\documentclass[ssy,preprint]{imsart} 

\RequirePackage{amsthm,amsmath}
\RequirePackage[numbers]{natbib}
\RequirePackage[colorlinks,citecolor=blue,urlcolor=blue]{hyperref}

\usepackage{amssymb,amsfonts,color,xcolor} 
\usepackage[english]{babel}

\usepackage[a4paper,bindingoffset=0.5cm,left=2cm,right=2cm,top=2.5cm,bottom=2cm,footskip=.8cm]{geometry} 

\arxiv{arXiv:1801.02999}

\startlocaldefs
\newcommand{\s}{^\star}
\newcommand{\vt}{\vartheta}
\newcommand{\tc}[1] {\textcolor{black}{#1}}
\newcommand{\tr}[1] {\textcolor{black}{#1}}

\newcommand{\tcp}[1]{\textcolor{black}{#1}}
\renewcommand{\bar}{\overline}
\newtheorem{lemma}{Lemma}
\newtheorem{corollary}{Corollary}
\newtheorem{assumption}{Assumption}
\newtheorem{theorem}{Theorem}
\newtheorem{remark}{Remark}
\newtheorem{example}{Example}

\endlocaldefs

\allowdisplaybreaks

\begin{document}

\vspace*{15pt} 

\renewcommand{\abstractname}{\enskip}
\setattribute{abstract}   {width} {.8\textwidth}

\begin{frontmatter}
\title{\vspace{3mm}\small Exact asymptotics for a multi-timescale model,\\
{with applications in modeling overdispersed customer streams}}
\runtitle{Exact asymptotics for a multi-timescale model}

\begin{aug}
\author{\fnms{Mariska} \snm{Heemskerk}
\thanksref{t1}
\ead[label=e1]{j.m.a.heemskerk@uva.nl}}
\& \,
\author{\fnms{Michel} \snm{Mandjes}
\ead[label=e2]{m.r.h.mandjes@uva.nl}} 



\affiliation{University of Amsterdam} 

\address{University of Amsterdam \\ Korteweg-de Vries Institute for Mathematics \\ Science Park 107 \\1098 XG  Amsterdam\\ The Netherlands.\\
\printead{e1}\\
\phantom{E-mail:\ }\printead*{e2}}
\end{aug}

\begin{abstract} 
{In} this paper we study the probability $\xi_n(u):={\mathbb P}\left(C_n\geqslant u n \right)$, with 
$C_n:=A(\psi_n B(\varphi_n))$ for L\'{e}vy processes $A(\cdot)$ and $B(\cdot)$, and $\varphi_n$ and $\psi_n$ non-negative sequences such that $\varphi_n \psi_n =n$ and $\varphi_n\to\infty$ as $n\to\infty$.  
Two timescale regimes are distinguished: a `fast' regime in which $\varphi_n$ is superlinear and a `slow' regime in which $\varphi_n$ is sublinear. 
We provide the exact asymptotics of $\xi_n(u)$ (as $n\to\infty$) for both regimes, relying on  change-of-measure arguments in combination with Edgeworth-type estimates. 
The asymptotics have an unconventional form: the exponent contains the commonly observed linear term, but may also contain sublinear terms (the number of which depends on the precise form of $\varphi_n$ and $\psi_n$).
To showcase the power of our results we include two examples, covering both the case where $C_n$ is lattice and non-lattice.
Finally we present numerical experiments that demonstrate the importance of taking into account the doubly stochastic nature of $C_n$ in a practical application related to customer streams in service systems;
they show that the asymptotic results obtained yield  highly accurate approximations, also in scenarios in which there is no pronounced timescale separation.
\end{abstract}



\end{frontmatter}

\section{Introduction, preliminaries, notation, literature}

Consider a scalar L\'evy process $A(\cdot)$, and (independently of $A(\cdot)$) an increasing scalar L\'evy process $B(\cdot)$. 
These are uniquely characterized by their {\it characteristic exponents}, which are defined as
\[\alpha(\vt) :=\log {\mathbb E}{\rm e}^{\vt A(1)},\:\:\:\beta(\vt) :=\log {\mathbb E}{\rm e}^{\vt B(1)},\]
and are in fact the logarithmic moment generating functions of $A(1)$ and $B(1)$.
Targeting at large-deviations asymptotics we impose the assumption that both characteristic exponents are finite in an open neighborhood of the origin, implying that all moments of $A(t)$ and $B(t)$ exist. 
Let  $\varphi_n$ and $\psi_n$ be non-negative sequences such that $\varphi_n \psi_n =n$ and $\varphi_n\to\infty$ as $n\to\infty$. 
Our objective is to find the {\it exact asymptotics} of
\[\xi_n(u):={\mathbb P}\left(A(\psi_n B(\varphi_n))\geqslant u n \right),\]
i.e., we wish to find a sequence $f_n$ such that $\xi_n(u)/f_n\to 1$ as $n\to\infty.$
We assume $u> ab$, with $a:={\mathbb E}\,A(1)>0$ and $b:={\mathbb E}\,B(1)>0$, such that the event under consideration becomes increasingly rare as $n\to\infty$. 
One special case has been studied in detail: the choice $\varphi_n=n$ and $\psi_n=1$ reduces $\xi_n(u)$ to ${\mathbb P}\left(A( B(n))\geqslant\, u n \right)$, whose exact asymptotics follow from \cite{BR} (using that $A(B(\cdot))$ is a L\'evy process). 
To the best of our knowledge, other cases have not been analyzed in the literature.

\subsubsection*{Effect of multiple timescales}
To get some intuition for the behavior of $\xi_n(u)$, we write it as
\[\xi_n(u) ={\mathbb P}\left(A\left(n \frac{B(\varphi_n)}{\varphi_n}\right)\geqslant u n \right).\] 
We proceed by explaining that there are two timescale regimes.
(i)~If $\varphi_n$ is superlinear, it is anticipated that $B(\varphi_n)/\varphi_n$ is close to $b$, such that  $\xi_n(u)$ resembles ${\mathbb P}(A(bn)\geqslant u n).$ 
We refer to this setting as the `fast regime', as the fluctuations of the process $B(\cdot)$ are so fast that it can be replaced by its mean value. 
(ii)~If on the contrary $\varphi_n$ is sublinear (which we will refer to as the `slow regime' for obvious reasons), 
then one may expect that the event of interest roughly looks like $an B(\varphi_n)/\varphi_n \geqslant u n$, and therefore $\xi_n(u)$ essentially behaves as ${\mathbb P}(a B(\varphi_n) \geqslant u\varphi_n)$.
The objective of this paper is to make these claims precise. 
Our main contribution is that we {succeed in identifying} the exact asymptotics of $\xi_n(u)$ as $n\to\infty$. 
These turn out to have a non-standard form, in the sense that  the exponent, in addition to the linear term that also appears in the classical asymptotics \cite{DZ}, may also contain sublinear terms.

To further investigate the two timescale regimes identified above, it is instructive to calculate the variance of $C_n:=A(\psi_n B(\varphi_n))$.
To this end, we first express the log-moment generating function (l-mgf) $\gamma_n(\cdot)$ of $C_n$ in terms of $\alpha(\cdot)$ and $\beta(\cdot).$
It requires a \tc{direct} computation to verify that
\[
\gamma_n(\vt)={\varphi_n } \,\beta\big( \alpha(\vt) \, \psi_n \big).
\]
Then it is direct that ${\mathbb E}\,C_n= n\alpha'(0)\beta'(0)= nab$ and
\begin{equation}
\label{VAR}{\mathbb V}{\rm ar}\,C_n =\gamma_n''(0) = n\psi_n \big(\alpha'(0)\big)^2\beta''(0)
+n\alpha''(0)\beta'(0)=n\psi_n \sigma^2_-
+n\sigma^2_+,\end{equation}
with  $\sigma^2_-:=a^2\beta''(0)$ and $\sigma^2_+:= \alpha''(0)b$.
In this decomposition of the variance, the {aforementioned regime} dichotomy is nicely reflected, as can be seen as follows. 
If $\varphi_n$ is superlinear (and hence $\psi_n$ vanishes), then the first term in the right-hand side of (\ref{VAR}) is small relative to the second term, and ${\mathbb V}{\rm ar}\,C_n$ {essentially} behaves as  $n\sigma^2_+$ (which is also the variance of $A(bn)$). 
On the other hand, if $\varphi_n$ grows sublinearly then so does $\psi_n$, so that in this case the first term of (\ref{VAR}) will dominate; as a consequence, ${\mathbb V}{\rm ar}\,C_n$ behaves as $n\psi_n \sigma^2_-$ (which equals the variance of $an B(\varphi_n)/\varphi_n$).

The above intuition can be translated in terms of  a central limit theorem by following a \tc{classical} approach. In the fast regime (in which $\varphi_n$ is superlinear), \tc{direct} computations reveal  that the characteristic function of 
\[D_n:= \frac{C_n-(ab)n}{\sqrt{n}\sigma_+}\]
converges to that of a standard Normal random variable. Hence, as a direct application of L\'evy's convergence theorem, $D_n$ converges in distribution to a standard Normal random variable. Likewise, in the slow regime (in which $\varphi_n$ is sublinear)
\[E_n:= \frac{C_n-(ab)n}{\sqrt{n\psi_n}\sigma_-}\]
converges in distribution to a standard Normal random variable. 
We conclude that the dichotomy {described} above manifests itself in a central limit context.

\subsubsection*{Large deviations, contributions}

In this paper {we assess} to what extent the above dichotomy carries over to a large-deviations setting. Based on the above reasoning, it is tempting to believe the following conjecture: in the fast regime $\xi_n(u)/{\mathbb P}(A(bn)\geqslant un)\to 1$ as $n\to\infty$, and in slow regime $\xi_n(u)/{\mathbb P}(a B(\varphi_n) \geqslant u\varphi_n)\to 1$ as $n\to\infty$. 
Our study, however, shows that this conjecture does {\it not} hold in general. 
In more detail, denoting by $k_n\sim\ell_n$ that $k_n/\ell_n\to 1$ as $n\to\infty$,
the main result is that we find explicit sequences $\lambda_{+,n}$ and $\lambda_{-,n}$ such that in the fast regime, as $n\to\infty$,
\[{\xi_n(u)}\sim{{\mathbb P}(A(bn)\geqslant un)\lambda_{+,n}},\]
and in the slow regime, as $n\to\infty$,
\[{\xi_n(u)}\sim{{\mathbb P}(a B(\varphi_n) \geqslant u\varphi_n)\lambda_{-,n}}.\]
Here, $\lambda_{+,n}$ and $\lambda_{-,n}$ do not necessarily equal 1, thus refuting the above conjecture. {Note that hereby the {exact asymptotics} of $\xi_n(u)$ are fully identified}, as sequences $\kappa_{+,n}$ and $\kappa_{-,n}$  such that ${\mathbb P}(A(bn)\geqslant un)\kappa_{+,n}\to 1$ and ${\mathbb P}(a B(\varphi_n) \geqslant u\varphi_n)\kappa_{-,n}\to1$ {are given in} \cite{BR}. (As an aside, we mention that in specific situations $\lambda_{+,n}$ and $\lambda_{-,n}$ {\it do}  equal 1, so that in those situations the conjecture applies; we return to this below.) 

{The resulting exact asymptotics of $\xi_n(u)$ can be written as} the product of a polynomial and an exponential part.
In the fast regime, the polynomial part is inversely proportional to $\sqrt{n}$, as was found before  in various related settings (such as the one studied in \cite{BR}). 
The exponential part, however, has a rather unusual shape: not only does the exponent contain a commonly observed term that is linear in $n$, in addition it consists of finitely many sublinear terms (the number of which depends on the specific form of $\varphi_n$ and $\psi_n$).
In the slow regime similar results apply, but with the role of $n$ taken over by $\varphi_n$, meaning that the polynomial part is inversely proportional to $\sqrt{\varphi_n}$ and the exponent is the sum of a term that is linear in $\varphi_n$ and finitely many terms that are $o(\varphi_n)$.

An immediate consequence of our result is that the conjecture we stated above {\it is} true if the timescales of both L\'evy processes are sufficiently separated; then there are no sublinear terms in the exponent (where `sublinear' means $o(n)$ in the fast regime and $o(\varphi_n)$ in the slow regime). More specifically, the conjecture holds in the fast regime if $n\psi_n\to 0$ (i.e., then $\lambda_{+,n}=1$), and in the slow regime if $\varphi_n/\psi_n\to 0$ (i.e., then $\lambda_{-,n}=1$). 
It is further remarked that {\it on a logarithmic scale} the conjecture always holds, albeit in the following (weaker) sense:
we show that in the fast regime
\[\lim_{n\to\infty}\frac{1}{n}\log{\xi_n(u)}=
\lim_{n\to\infty}\frac{1}{n}\log
{\mathbb P}(A(bn)\geqslant un),\]
whereas in the slow regime
\[\lim_{n\to\infty}\frac{1}{\varphi_n}\log{\xi_n(u)}=
\lim_{n\to\infty}\frac{1}{\varphi_n}\log
{\mathbb P}(a B(\varphi_n) \geqslant u\varphi_n).\]

\begin{example} \label{exf}{\em
The leading example, which we will visit several times in this paper,
is that of $\varphi_n \sim n^f$ and $\psi_n \sim n^{1-f}$ for some $f>0$. 
The fast regime corresponds to $f>1$ and the slow regime to $f\in(0,1)$. 
The above criterion yields that $\lambda_{+,n}=1$ if $f>2$ and $\lambda_{-,n}=1$ if $f\in(0,\frac{1}{2}).$ \hfill$\Diamond$}
\end{example}

\begin{remark} \label{phin}\tc {\em We mentioned before that the case $\varphi_n = n$ can be dealt with in a straightforward way. For the sake of completeness, we include this argumentation here.  In this case $\gamma_n(\vt) = n \beta(\alpha(\vt))$, which effectively means that $C_n$ can be written as the sum of $n$ i.i.d.\ random variables, bringing us directly in the setting of \cite{BR}. With $\vt\s$ solving $\beta'(\alpha(\vt))\alpha'(\vt)=u$, one obtains 
\[\xi_n(u)\sim\frac{1}{\vt\s\sigma_0 \sqrt{2\pi n}} \exp\left((\beta(\alpha(\vt\s))-\vt\s u)n\right),\]
where $\sigma_0:= \beta''(\vt\s)(\alpha'(\vt\s))^2+ \beta'(\alpha(\vt\s))\alpha''(\vt\s).$
\hfill$\Diamond$}
\end{remark}

We proceed with a few words on the approach followed. 
In the large-deviations analysis a crucial role is played by the `twisting factor', i.e.,
the solution $\vt = \vt_n$ of the equation $un =  \gamma_n'(\vt)$, or
\begin{equation}\label{twist}u= {\beta'\big(\alpha(\vt)\,\psi_n \big)}\,{\alpha'(\vt)}.\end{equation}
As $u >ab$, it follows that $\vt_n$ is positive. 
In the fast regime, $\vt_n$ is close to the solution $\vt\s$ of $u= b\alpha'(\vt\s)$; in the slow regime $\vt_n$ resembles $\tau\s/\psi_n$ where $\tau\s$ solves $u=\beta'(a\tau\s)a$. 
Our proofs rely on a change of measure via an exponential twist that is based on the solution of (\ref{twist}). 
By and large, the proof presented in \cite{DZ} underlying the Bahadur-Rao \cite{BR} result is followed: 
{first we capture the exponential part of $\xi_n(u)$ by applying a change of measure},
{after which the polynomial part of $\xi_n(u)$ is identified using delicate Berry-Ess\'een-based calculations (that are considerably more involved than the ones needed in the `classical' Bahadur-Rao case).}
In line with the proof in \cite{DZ}, the case where \tc{the underlying L\'evy processes} are lattice has to be dealt with separately.

\tcp{Concluding the technical part of the paper, we want to apply the gained insights in a practical example, hence demonstrating its use in OR. In this part we (i)~provide more background on the rationale behind modeling overdispersion, (ii)~we point out how the asymptotics identified in our paper can be translated into approximations (for an unscaled model), (iii)~we assess the accuracy of such approximations.}

\subsubsection*{Literature and motivation}
Our work fits in the tradition of large-deviations asymptotics of sample-mean related quantities in the light-tailed setting. Without attempting to provide an exhaustive overview, we give a number of key references. 

In the classical framework, $S_n$ is defined as the sum of i.i.d.\ random variables $X_1, \dots, X_n$, where the $X_i$ are assumed to have a finite moment generating function in a neighborhood of the origin. 
{The exceedance probability $\chi_n(u):={\mathbb P}(S_n\geqslant un)$ is the object of study, with $u> {\mathbb E}X_1$.}
Cram\'er~\cite{CRA} characterized a function $I(u)$ such that $\tfrac{1}{n}\log \chi_n(u)\to -I(u)$ as $n\to\infty$; this type of results is commonly referred to as {\it logarithmic asymptotics}.  
The  proof technique used relied on  change-of-measure argumentation that has been applied extensively since then. 

Cram\'er's seminal work was extended in several directions.  In \cite{CHE} a uniform upper bound on $\chi_n(u)$, generally known as the {\it Chernoff bound},  was derived. 
We refer to 
\cite{BZ} for a generalization of Cram\'er's logarithmic asymptotics to infinite-dimensional topological vector spaces.
Cram\'{e}r's theorem was also extended to sums of dependent random variables \cite{PLA} and vectors \cite{ELL,GAR}.

Whereas the above results focus on logarithmic asymptotics, another strand of research addresses exact asymptotics. 
In this respect we mention the pioneering paper \cite{BR} that shows that, under the conditions of Cram\'er's result, $\chi_n(u) \,{\rm e}^{nI(u)}\sqrt{n}$ converges to a positive constant as $n\to\infty$; 
in \cite{BH} a similar result had already been derived for the case that the $X_i$ are lattice. 
We refer to e.g.\ \cite{CS} for exact asymptotics in the vector-valued case, and to \cite{HOG} for a uniform framework covering both the CLT regime and large deviations. 

Our investigations were motivated by models recently suggested {in order} to incorporate {\it overdispersion}. 
The reason for developing such models was the observation that in various types of service systems \cite{BAS,WHI} the customer arrival process is intrinsically more variable than the traditionally used Poisson process. 
An approach to {overcome} the lack of overdispersion was proposed in \cite{HLM}: extra variability is produced by periodically {\it resampling} the Poisson arrival rate. As a consequence, when resampling every unit of time,
the number of arrivals in $[0,m]$ denoted by \tc{$C(m)$} (for $m\in{\mathbb N}$) is Poisson distributed with (random) parameter 
$\sum_{i=1}^m \Lambda_i$,
where the $\Lambda_i$ are i.i.d.\ non-negative random variables. 
Assuming for simplicity that the $\Lambda_i$ are  integer-valued as well, this means that
\begin{equation}\tc{C(m)} \stackrel{\rm d}{=} \sum_{j=1}^{\sum_{i=1}^m \Lambda_i} X_j,\label{TTP}\end{equation}
with $X_j$ a sequence of i.i.d.\ Poisson random variables with parameter 1. We thus have introduced overdispersion, as desired; more specifically, we have 
\[{\mathbb V}{\rm ar}\,C(m) =m\,{\mathbb E}\Lambda\,{\mathbb V}{\rm ar}X+m\, ({\mathbb E}X)^2{\mathbb V}{\rm ar}\Lambda\geqslant  m\,{\mathbb E}\Lambda\,{\mathbb V}{\rm ar}X= m\,{\mathbb E}\Lambda ={\mathbb E}\,C(m).\]
making the process more variable than the Poisson process.
\tcp{Observe} that the `two-timescale random walk' (\ref{TTP}) we introduced can be seen as the discrete counterpart of the two-timescale L\'evy processes considered in the present paper. In \cite{HLM} the above two-timescale random walk model is studied under certain scalings (comparable to the scalings with $\varphi_n$ and $\psi_n$ that are imposed in the present paper). 
The results in \cite{HLM} include 
logarithmic asymptotics, which are for special cases refined to exact asymptotics in~\cite{HKM}.


\subsubsection*{Organization} 
The fast regime (in which $\varphi_n$ grows superlinearly) is covered by Section \ref{FR}, followed by the slow regime (in which $\varphi_n$ grows sublinearly) in Section \ref{SR}. Examples are presented in Section~\ref{EX}. Section \ref{OVERD} describes how our findings can be applied when modeling overdispersed customers streams \tcp{and finally, Section \ref{DISC} puts our results into perspective and identifies directions for follow-up research.}

\section{Fast regime} \label{FR}
In this section we analyze the case that $\varphi_n$ is superlinear, entailing that $\psi_n \to0$ as $n\to\infty$; this corresponds to $f>1$ in the context of Example \ref{exf}.
\tc{The approach followed echoes the proof of e.g.\ \cite[Thm.\ 3.7.4]{DZ}, and comprises the following elements:
\begin{itemize}
\item[$\circ$]
We first identify in Section 2.1 the `twisting factor'  $\vt_n$ (i.e., the solution of (\ref{twist})). It means that twisting $C_n$ by $\vt_n$ leads to a random variable with mean $un$. More concretely, we define a measure ${\mathbb Q}_n$ through (in self-evident notation)
\[{\mathbb Q}_n(C_n\in{\rm d}x) = {\mathbb P}(C_n\in{\rm d}x) \tcp{\frac{\exp(\vt_n x)}{\exp \gamma_n(\vt_n)}},\]
where it holds that  (in self-evident notation)
\[{\mathbb E}_{{\mathbb Q}_n} C_n=\int_{-\infty}^\infty x\, {\mathbb Q}_n(C_n\in{\rm d}x) = un.\]
\item[$\circ$]
Then we rewrite in Section 2.2 the probability $\xi_n(u)$ using the $\vt_n$-twisted version of $C_n$. As in \cite[Ch.\ XIII]{AS}, we can express $\xi_n(u)$ into a mean under ${\mathbb Q}_n$:
\begin{equation}\label{COMM}\xi_n(u) ={\mathbb P}(C_n\geqslant un) = {\mathbb E}_{{\mathbb Q}_n}(L(C_n)1\{C_n\geqslant un\}),\end{equation}
for an appropriately chosen $L(\cdot)$ (which can be interpreted as a likelihood ratio). 
\item[$\circ$] The right-hand side of (\ref{COMM}) turns out to consist of an exponential part and a polynomial part  (in $n$, that is). 
Section 2.3 analyzes the exponential part. As opposed to the standard Bahadur-Rao result \cite[Thm.\ 3.7.4]{DZ}, the exponent potentially also contains sublinear terms. 
\item[$\circ$] The polynomial part in the right-hand side of (\ref{COMM}) is technically the most demanding element of the proof, and relies on Berry-Ess\'een-based arguments; see Section 2.4.
\item[$\circ$]  In Section 2.5 all elements are put together, and our result is stated. 
\end{itemize}}
The formal assumption we impose on $\psi_n$ is the following.
\begin{assumption} \label{ass1} 
The sequence $\psi_n$ satisfies
\[\limsup_{n\to\infty} \frac{\log \psi_n}{\log n}<0.\]
\end{assumption}
This assumption means that there is an $\varepsilon>0$ such that $\psi_n<n^{-\varepsilon}$, and hence $\varphi_n> n^{1+\varepsilon}.$ We observe that $\varphi_n$ is superlinear. 

\tc{In the main theorem of this section, i.e., Thm.\ \ref{THM1}, we assume that $A(\cdot)$ is non-lattice.
To establish the result for the case where $A(\cdot)$ is lattice (covering e.g.\ Poisson processes), a minor adaptation needs to be made. More specifically, the steps of Sections 2.1 up to 2.3 apply to the lattice as well as non-lattice case, whereas in Section~2.4 in the lattice case slightly different bounds have to be used (fully analogously to the treatment of the lattice and non-lattice cases in the proof of \cite[Thm. 3.7.4]{DZ}).
The result for the lattice case is discussed separately in Remark \ref{REM1}. We refer to Sections \ref{EX} and \ref{OVERD} for illustrative examples covering both the lattice and non-lattice case.}

\subsection{Analysis of twisting factor} \tc{As $\psi_n\to 0$ as $n\to\infty$}, the twisting factor
$\vt_n$ converges to the solution $\vt\s$ of $b \alpha'(\vt)=u$,
where $b={\mathbb E}B(1)=\beta'(0)$.
\tc{To establish the exact asymptotics of $\xi_n(u)$, it turns out that we need the an expansion of $\vt_n$; we start by arguing that $\vt_n$ has the form
\begin{equation}
\label{thetan}\vt_n = \vt\s +\sum_{k=1}^\infty v_k  \psi_n ^{\:k},
\end{equation}
and then we point out how to determine the coefficients $v_k$.}

\begin{itemize}
\item[$\circ$] \tc{The main idea is that one can (implicitly) define a function $\vt_+(x)$ as the solution of the equation
\[\beta'(\alpha(\vt)x)\alpha'(\vt) =u;\]
such a solution is unique due to the fact that the left-hand side of the previous display is the derivative of a convex function (hence increasing). \tcp{Observe} that $\vt_n$ (solving $\gamma_n'(\vt_n)=u$) equals $\vt_+(\psi_n)$; recall that $\psi_n\to 0$ as $n\to\infty.$}

\noindent
\tc{Observe that $\vt_+(0)=\vt\s$ (with $\vt\s$, as defined above, solving $b\alpha'(\vt\s)=u$). 
\tcp{Writing} down the Taylor expansion (at $x=0$) of the implicitly defined function $\vt_+(x)$:
\[\vt_+(x) =\sum_{k=0}^\infty v_k x^k,\]
with $v_0=\vt\s$,
the $v_k$ can be determined in terms of the derivatives of $\alpha(\cdot)$ and $\beta(\cdot)$
(which are well-defined). Below we provide a constructive procedure that identifies all these coefficients. Upon combining the above elements, we thus obtain the expansion (\ref{thetan}):
\[\vt_n=\vt_+(\psi_n)= \sum_{k=0}^\infty v_k \psi_n^{\,k}.\]}

\item[$\circ$]
We present a constructive procedure that yields the coefficients $v_k$.
To this end, first observe that, by evaluating $\alpha(\cdot)$ as a Taylor series {around $\vt\s$},
\[\beta'\left(\sum_{\ell=0}^\infty \frac{\alpha^{(\ell)}(\vt\s)}{\ell!}\left(\sum_{k=1}^\infty v_k  \psi_n ^{\:k}\right)^\ell \psi_n\right)\cdot\sum_{\ell=0}^\infty \frac{\alpha^{(\ell+1)}(\vt\s)}{\ell!}\left(\sum_{k=1}^\infty v_k  \psi_n ^{\:k}\right)^\ell=u ,\]
where, by evaluating $\beta(\cdot)$ as a Taylor series {around $0$},
\[\beta'(\vt) = \sum_{m=0}^\infty \frac{\beta^{(m+1)}(0)}{m!} \vt^m.\]
This equation effectively determines all $v_k$, \tcp{which can be found by} grouping together the appropriate terms \tc{(i.e., terms with the same power). } \tc{We demonstrate this procedure for $v_1$. }
\tc{Up to terms that are $o(\psi_n)$,
\[\left(\beta'(0) + {\beta''(0)}\alpha(\vt\s)\psi_n\right)
\left(\alpha'(\vt\s)+ {\alpha''(\vt\s)}v_1\psi_n\right)=u.\]
\tcp{Using that} $\beta'(0)\alpha'(\vt\s)=b\alpha'(\vt\s)=u$, we obtain}
\[\tc{v_1 = -\frac{\alpha(\vt\s)\,\alpha'(\vt\s)}{\alpha''(\vt\s)} \frac{\beta''(0)}{\beta'(0)}.}\]
The $v_i$ for $i\in\{2,3\ldots\}$ can be computed in the exact same way, but this does not lead to clean expressions; for $k=2$ we obtain
\[
v_2= -\frac12\frac{\beta^{(1)}(0)\alpha^{(3)}(\theta^*)v_1^2 +\beta^{(2)}(0)(\alpha(\theta^*)\alpha^{(2)}(\theta^*) + (\alpha^{(1)}(\theta^*))^2)v_1}{\beta^{(1)}(0) \alpha^{(2)}(\theta^*)}.
\]
\end{itemize}

\subsection{Change of measure}
The next step is to rewrite the probability of interest, relying on the usual change-of-measure procedure. 
This effectively means that we let ${\mathbb Q}_n$ correspond to twisting {the distribution of} $C_n$ by the solution $\vt_n$ of~(\ref{twist}). 
\tc{The l-mgf of $C_n$ under ${\mathbb Q}_n$ can thus be expressed in terms 
of the l-mgf of $C_n$ under the original measure:}
\[\gamma^{\mathbb Q_n}_n(\vt) := \gamma_n(\vt+\vt_n)-\gamma_n(\vt_n).\]
Later on we need the mean and variance of $C_n$ under ${\mathbb Q}_n$. 
\tc{From the definition of $\vt_n$, it follows that ${\mathbb E}_{{\mathbb Q}_n}\,C_n = u n$. }
Differentiating once more, we obtain
\[{\mathbb V}{\rm ar}_{{\mathbb Q}_n}C_n = n\psi_n \beta''\big(\alpha(\vt_n)\psi_n\big)\big(\alpha'(\vt_n)\big)^2 + n \beta'\big(\alpha(\vt_n)\psi_n\big)\alpha''(\vt_n).\]
 
The next step is to use the change of measure to rewrite $\xi_n(u)$.
By applying the definition of ${\mathbb Q}_n$, we find the identity
\[\xi_n(u)= {\mathbb E}_{{\mathbb Q}_n}\left( {\rm e}^{ \gamma_n(\vt_n)-\vt_nC_n}\,1\{C_n \geqslant u n \}\right),\]
realizing that ${\rm e}^{ \gamma_n(\vt_n)-\vt_n C_n}$ {plays the role of} the likelihood ratio ${\rm d}{\mathbb P}/{\rm d}{{\mathbb Q}_n}$. 
\tc{On the event $C_n\geqslant un$, $C_n$ will typically be relatively close to $un$ under the new measure ${\mathbb Q}_n$ (recall ${\mathbb E}_{{\mathbb Q}_n}\,C_n = u n$).}
To exploit this, we define, with $\sigma_+^{\mathbb Q}:=\sqrt{b\alpha''(\vt\s)}$,
\[\bar D_n := \frac{C_n -un}{\sqrt{n}\sigma_+^{\mathbb Q}},\]
which is a random variable that has, by construction, mean $0$ under ${\mathbb Q}_n$. \tc{(Due to the $\sqrt{n}$ in the denominator, it is anticipated that, under ${\mathbb Q}_n$, $\bar D_n$ can be approximated by a standard Normal random variable; we come back to this idea below.)}
We thus obtain that, for all $n$,
\begin{equation}\label{ident_fast}
\xi_n(u) = {{\rm e}^{\gamma_n(\vt_n)-\vt_n u n} \, \Delta_n},\:\:\:\mbox{with}\:\:\Delta_n:= {\mathbb E}_{{\mathbb Q}_n}\left( {\rm e}^{-\vt_n\sigma_+^{\mathbb Q}\sqrt{n}\, \bar D_n}\,1\{\bar D_n \geqslant 0 \}\right).\end{equation}
Consequently, we are left with analyzing the exponential factor $\delta_n:=\exp(\gamma_n(\vt_n)-\vt_n u n)$ and the expectation $\Delta_n$, as $n$ grows large.

\subsection{Analysis of $\delta_n$ as $n\to\infty$}
We proceed by analyzing the exponent in the expansion~(\ref{ident_fast}), i.e., $\gamma_n(\vt_n)-\vt_n u n$.
Define 
\[
m_+:=\sup\left\{k\in{\mathbb N}: \liminf_{n\to\infty}\varphi_n\psi_n^{\,k}>0\right\}; 
\]
note that $m_+\geqslant 1$.
Our claim is the following.
\begin{lemma}\label{EXP}As $n\to\infty$, for constants ${\bar v}_k$,
\[
\gamma_n(\vt_n)-\vt_n u n=
\big(\,b\alpha(\vt\s) -\vt\s u\,\big)n +\sum_{k=2}^{m_+} {\bar v}_k \varphi_n\psi_n^{\, k}+o(1),\]
where the empty sum is defined as $0$.
\end{lemma}
The validity of this claim (as well as the procedure to determine the coefficients ${\bar v}_k$)
can be demonstrated as follows.
Relying on Taylor expansions, we obtain
\begin{align*}
\gamma_n(\vt_n)-\vt_n  un&= \varphi_n\,\beta\left(\sum_{\ell=0}^\infty \frac{\alpha^{(\ell)}(\vt\s)}{\ell!}\left(\sum_{k=1}^\infty v_k  \psi_n ^{\:k}\right)^\ell \psi_n\right)-\left(\vt\s +\sum_{k=1}^\infty v_k  \psi_n ^{\:k}\right)un.
\end{align*}
The claim for $m_+=1$ can be directly verified.
Now consider the case $m_+=2$; any higher value of $m$  can be dealt with fully analogously. 
For $m_+=2$, the sequence ${ \varphi_n\,\psi_n ^{\:k}}$ converges to $0$ when $k> 2$, whereas $\varphi_n\,\psi_n^{\, 2}=n\psi_n$ {stays away from} $0$;
think of e.g.\ $\psi_n=n^{-2/3}$.
Hence we obtain that, as $n\to\infty$,
\begin{align*} 
\gamma_n(\vt_n)-\vt_n u n&= \left(b\alpha(\vt\s)-\vt\s u\right)n +  v_1\left(\frac{\beta''(0)}{2}\alpha(\vt\s)+b\alpha'(\vt\s)-u\right) \varphi_n\psi_n^{\, 2}+o(1).
\end{align*}
It is clear that for $m_+=3$, an additional term needs to be included. 
In general, using this procedure any value of $m_+$ can be dealt with.

\subsection{Analysis of $\Delta_n$ as $n\to\infty$}
We are left with analyzing $\Delta_n$ for $n$ large. Our objective is to prove that $\sqrt{n}\Delta_n$ converges to a positive constant as $n\to\infty.$ 
Mimicking the line of reasoning of \cite[Eqn. (3.7.7)]{DZ}, {we apply} integration by parts:
\begin{align}\nonumber
\Delta_n = \int_0^\infty {\rm e}^{-\vt_n\sigma_+^{\mathbb Q}\sqrt{n} x}{\mathbb Q}_n(\bar D_n\in{\rm d}x)
&=\sqrt{n}\vt_n \sigma_+^{\mathbb Q} \int_0^\infty {\rm e}^{-\vt_n\sigma_+^{\mathbb Q}\sqrt{n} x} \big({\mathbb Q}_n(\bar D_n\leqslant x)-
{\mathbb Q}_n(\bar D_n\leqslant 0)\big){\rm d}x \\
&=\vt_n\sigma_+^{\mathbb Q}  \int_0^\infty {\rm e}^{-\vt_n\sigma_+^{\mathbb Q} x} \big({\mathbb Q}_n(\bar D_n\leqslant x/\sqrt{n})-
{\mathbb Q}_n(\bar D_n\leqslant 0)\big){\rm d}x.
\label{alte}
\end{align}
As we {will} intensively rely on this representation of $\Delta_n$, 
an important role in our argumentation is played by the probability distribution
\[{\mathbb Q}_n(\bar D_n\leqslant x) = {\mathbb Q}_n\left(\frac{A(\psi_nB(\varphi_n))-un}{\sqrt{n}\sigma_+^{\mathbb Q}}\leqslant x\right),\]
where it is noted that ${\mathbb E}_{{\mathbb Q}_n} A(\psi_nB(\varphi_n)) = u n$.
Estimates for this distribution can be established, essentially as variants of the classical Edgeworth expansion. As such expansions follow by applying \tc{well-established} techniques, we restrict ourselves to providing the main steps of the derivation in \tcp{Appendix \ref{Edge}}. An excellent introduction on the Edgeworth expansion is provided in \cite[Ch. II]{HALL}; see also \cite[Ch. IV]{BDC}.

$\circ$\:\:In the case where $\lim_{n \to \infty }\psi_n \sqrt{n} = 0$ (which corresponds to $f>\frac{3}{2}$ in the context of Example \ref{exf}), we have, as pointed out in \tcp{Appendix \ref{Edge}}, as $n\to\infty$, 
\begin{equation}
\label{B1}\sqrt{n}\sup_x\left({\mathbb Q}_n(\bar D_n\leqslant x)- \Phi(x)+\phi(x) H_2(x)\frac{\kappa_+}{\sqrt{n}}\right)\to 0,\:\:\:\kappa_+:=\frac{1}{6}\frac{b\alpha'''(\vt\s)}{(\sigma_+^{\mathbb Q})^3};\end{equation}
cf.\ Eqn.\ (\ref{bound1fast}), where $H_2(x)=x^2-1$, 
$\phi(\cdot)$ denotes the probability density function of a standard Normal random variable, and $\Phi(\cdot)$ the corresponding cumulative distribution function. 
We are now in a position to prove that $\sqrt{n}\Delta_n$ converges to a constant; we provide the upper bound, but the lower bound follows fully analogously. By virtue of (\ref{B1}), for any $\varepsilon>0$ and $n$ sufficiently large, 
\begin{align*}
\sqrt{n}\Delta_n \leqslant &\: \sqrt{n}\vt_n \sigma_+^{\mathbb Q} \int_0^\infty {\rm e}^{-\vt_n\sigma_+^{\mathbb Q} x} \left(\Phi(\frac{x}{\sqrt{n}})-\Phi(0)\right)
{\rm d}x\,-\\
&\:
\sqrt{n}\vt_n\sigma_+^{\mathbb Q}  \int_0^\infty {\rm e}^{-\vt_n\sigma_+^{\mathbb Q} x} \left(\phi(\frac{x}{\sqrt{n}})H_2(\frac{x}{\sqrt{n}})-\phi(0)H_2(0)\right)\frac{\kappa_+}{\sqrt{n}}
{\rm d}x\,
+\,\varepsilon\vt_n \sigma_+^{\mathbb Q} \int_0^\infty {\rm e}^{-\vt_n\sigma_+^{\mathbb Q} x} 
{\rm d}x.
\end{align*}
Let us start by evaluating the first term on the right-hand side. 
It can be rewritten as
\begin{align*}
\sqrt{n}\vt_n\sigma_+^{\mathbb Q}  \int_0^\infty {\rm e}^{-\vt_n\sigma_+^{\mathbb Q} x} \int_0^{x/\sqrt{n}} \frac{1}{\sqrt{2\pi}}{\rm e}^{-y^2/2}{\rm d}y\,
{\rm d}x
=&\,\sqrt{n}\vt_n \sigma_+^{\mathbb Q} \int_0^\infty\frac{1}{\sqrt{2\pi}}{\rm e}^{-y^2/2} \int_{y\sqrt{n}}^\infty {\rm e}^{-\vt_n\sigma_+^{\mathbb Q} x} {\rm d}x\,
{\rm d}y\\=&\, {\sqrt{n}} \int_0^\infty\frac{1}{\sqrt{2\pi}}{\rm e}^{-y^2/2}  {\rm e}^{-\vt_n\sigma_+^{\mathbb Q} \sqrt{n} \,y} 
{\rm d}y\\
=&\, {\sqrt{n}}
\exp\left(\frac{1}{2}( \vt_n\sigma_+^{\mathbb Q})^2   n\right) \left(1-\Phi(\vt_n\sigma_+^{\mathbb Q} \sqrt{n})\right).
\end{align*}
Using the known limit $x(1-\Phi(x))/\phi(x)\to 1$ as $x\to\infty$, and $\vt_n\to\vt\s$, we have that the expression in the previous display converges to
\begin{equation}\label{constantefast}
\big({\vt\s \sigma_+^{\mathbb Q} \sqrt{2\pi}}\big)^{-1}.
\end{equation}
Using that $H_2(x)= x^2-1$ the second term can be written as the sum of
\begin{align*}
t_{+,n}^{(1)}&:=\,\sqrt{n}\vt_n \sigma_+^{\mathbb Q} \int_0^\infty {\rm e}^{-\vt_n\sigma_+^{\mathbb Q} x} \left(\phi(0)-\phi(\frac{x}{\sqrt{n}})\right)\frac{\kappa_+}{\sqrt{n}}
{\rm d}x,\\
t_{+,n}^{(2)}&:=\,
\sqrt{n}\vt_n \sigma_+^{\mathbb Q} \int_0^\infty {\rm e}^{-\vt_n\sigma_+^{\mathbb Q} x} \phi(\frac{x}{\sqrt{n}})\frac{x^2}{n}\frac{\kappa_+}{\sqrt{n}}
{\rm d}x.\end{align*}
From (i)~$\sqrt{n}\,(\phi(0)-\phi({x}/{\sqrt{n}}))\to x\phi'(0)=0$, (ii)~$\vt_n\to\vt\s$ as $n\to\infty$, and (iii)~the fact that $\phi'(\cdot)$ is bounded, applying the dominated convergence theorem yields that $t_{+,n}^{(1)}\to 0.$ Due to (i)~$\phi(x/\sqrt{n})\leqslant 1/\sqrt{2\pi}$, (ii)~$\vt_n\to\vt\s$ as $n\to\infty$, and (iii)
\[\limsup_{n\to\infty}\int_0^\infty  {\rm e}^{-\vt_n\sigma_+^{\mathbb Q} x} x^2\,{\rm d}x = \limsup_{n\to\infty}\frac{2}{(\vt_n\sigma_+^{\mathbb Q})^3}<\infty,\]
we also have $t_{+,n}^{(2)}\to 0.$

The third term equals $\varepsilon$, which can be made arbitrarily small. As mentioned before, \tc{by the same token it can be  verified that the corresponding lower bound applies as well}. We thus conclude that $\sqrt{n}\Delta_n$ converges to the constant in (\ref{constantefast}).

$\circ$\:\:We proceed with the case  $\liminf_{n\to \infty} \psi_n \sqrt{n} > 0$  (which simplifies  to $f\in(1,\frac{3}{2}]$ in the context of Example \ref{exf}). 
From Eqn.\ (\ref{bound2fast}) in \tcp{Appendix \ref{Edge}} we have
\begin{equation}
\label{B2}\sqrt{n}\sup_x\left({\mathbb Q}_n(\bar D_n\leqslant x)- \Phi(x)+\phi(x) \left(H_1(x)\sum_{k=1}^{m_+}c_k\psi_n^{\, k} +H_2(x)\frac{\kappa_+}{\sqrt{n}}\right)\right)\to 0.\end{equation}
By (\ref{B2}), for any $\varepsilon>0$ and $n$ sufficiently large, using that $H_1(x)=x$ and hence $H_1(0)=0$,
\begin{align*}
\sqrt{n}\Delta_n \leqslant &\:\sqrt{n}\vt_n \sigma_+^{\mathbb Q} \int_0^\infty {\rm e}^{-\vt_n\sigma_+^{\mathbb Q} x} \left(\Phi(\frac{x}{\sqrt{n}})-\Phi(0)\right)
{\rm d}x\,{-}\\
&\:\sqrt{n}\vt_n \sigma_+^{\mathbb Q} \int_0^\infty {\rm e}^{-\vt_n\sigma_+^{\mathbb Q} x}  \phi(\frac{x}{\sqrt{n}})\frac{x}{\sqrt{n}}\left(\sum_{k=1}^{m_+}c_k\psi_n^k \right) 
{\rm d}x\,{-}\,\\
&\:\sqrt{n}\vt_n \sigma_+^{\mathbb Q} \int_0^\infty {\rm e}^{-\vt_n\sigma_+^{\mathbb Q} x} \left(\phi(\frac{x}{\sqrt{n}})H_2(\frac{x}{\sqrt{n}})-\phi(0)H_2(0)\right)\frac{\kappa_+}{\sqrt{n}}
{\rm d}x\,+\,\varepsilon\vt_n \sigma_+^{\mathbb Q} \int_0^\infty {\rm e}^{-\vt_n\sigma_+^{\mathbb Q} x} 
{\rm d}x
\end{align*}
The first, third, and fourth term can be dealt with as in the case $\psi_n \sqrt{n}\to 0$. 
Due to (i)~$\phi(x)\leqslant 1/\sqrt{2\pi}$ for all $x$,  (ii)~$\vt_n\to\vt\s$ as $n\to\infty$, (iii)~$\psi_n^k\to 0$ as $n\to\infty$ for any $k\in\{1,\ldots,m_+\}$, and (iv)
\[\limsup_{n\to\infty}\int_0^\infty  {\rm e}^{-\vt_n\sigma_+^{\mathbb Q} x} x\,{\rm d}x = \limsup_{n\to\infty}\frac{1}{(\vt_n\sigma_+^{\mathbb Q})^2}<\infty,\]
the second term vanishes. 
It follows that again $\sqrt{n}\Delta_n$ converges to the constant in (\ref{constantefast}). 
We have thus established the following  lemma.

\begin{lemma} \label{lemD} As $n\to\infty$,
\[\sqrt{n}\Delta_n \to \big({\vt\s \sigma_+^{\mathbb Q} \sqrt{2\pi}}\big)^{-1}.\]
\end{lemma}

\subsection{Result}

Upon combining Lemmas \ref{EXP} and \ref{lemD}, as presented in the previous subsections, we have arrived at the following result.

\begin{theorem}\label{THM1}
As $n\to\infty$, under Assumption \ref{ass1}, \tc{for non-lattice $A(\cdot)$,}
\[\xi_n(u)\sim \frac{1}{\vt\s \sigma_+^{\mathbb Q} \sqrt{2\pi n}}
\exp\left(\big(\,b\alpha(\vt\s) -\vt\s u\,\big)n +\sum_{k=2}^{m_+} {\bar v}_k \varphi_n\psi_n^{\,k}\right).\]
\end{theorem}

An immediate consequence of this theorem is that  $\xi_n(u)$ behaves as ${\mathbb P}(A(bn)\geqslant un)$
when $\varphi_n\psi_n^2 = n\psi_n\to 0$. This is for instance the case when $\psi_n=n^{\zeta}$ for $\zeta<-1$, or, in the setting of Example~\ref{exf}, $f>2$. It reflects that the timescale of 
$B(\cdot)$ is so much faster than that of $A(\cdot)$, that it can be replaced by its mean. In addition, it implies that
the rough (logarithmic) asymptotics are not affected by the choice of $\psi_n$ (as long as Assumption \ref{ass1} is fulfilled). 
These findings are summarized in the following corollary.

\begin{corollary}\label{C1}
If  $\varphi_n\psi_n^2 = n\psi_n\to 0$ as $n\to\infty$, then, under Assumption \ref{ass1}, \tc{for non-lattice $A(\cdot)$,}
\[\xi_n(u) \sim {\mathbb P}(A(bn)\geqslant un) \sim  \frac{1}{\vt\s \sigma_+^{\mathbb Q} \sqrt{2\pi n}}
\exp\left(\big(\,b\alpha(\vt\s) -\vt\s u\,\big)n \right).\]
As $n\to\infty$, under Assumption \ref{ass1}, \tc{for non-lattice $A(\cdot)$,}
\[\frac{1}{n}\log\xi_n(u) \to b\alpha(\vt\s) -\vt\s u.\]
\end{corollary}

\begin{remark}{\em 
\label{REM1} \tc{So far we throughout assumed that the L\'evy processes $A(\cdot)$ be non-lattice, thus ruling out e.g.\ Poisson processes. 
With a minor adaptation, however, the lattice case can be dealt with as well.
Let, for some $x_0$ and $d$ and any $t\geqslant 0$, the random variable $d^{-1}(A(t) -x_0)$  be integer almost surely, and let $d$ be the largest number with this property; e.g.\ in the Poisson process case $d=1$. 
Then, following the proof in e.g.\ \cite[Thm.\ II.7.4.b]{DZ}, under Assumption \ref{ass1}, we find the following counterpart of Thm.\ \ref{THM1}:
\[\xi_n(u)\sim \frac{d}{1-{\rm e}^{-\vt\s d}}\frac{1}{ \sigma_+^{\mathbb Q} \sqrt{2\pi n}}
\exp\left(\big(\,b\alpha(\vt\s) -\vt\s u\,\big)n +\sum_{k=2}^{m_+} {\bar v}_k \varphi_n\psi_n^{\,k}\right),\]
as $n\to\infty$. This behavior is consistent with that of Thm.\ \ref{THM1} when taking $d\downarrow 0.$
\hfill$\Diamond$}}
\end{remark}

\section{Slow regime}\label{SR}
In this section we analyze the case that $\varphi_n$ grows sublinearly, implying that $\psi_n \to\infty$ as $n\to\infty$ (also sublinearly). 
As before, we first characterize the solution $\vt_n$ of (\ref{twist}), then we rewrite $\xi_n(u)$ using the $\vt_n$-twisted version of $C_n$, and finally we determine the corresponding exact asymptotics.
The formal assumption we impose on $\psi_n$ in this section is the following.
\begin{assumption}\label{ass2} 
The sequence $\psi_n$ satisfies
\[
0 < \liminf_{n\to\infty} \frac{\log \psi_n}{\log n} \leqslant \limsup_{n\to\infty} \frac{\log \psi_n}{\log n} <1.
\]
\end{assumption}
The first inequality of this assumption {ensures} that there is an $\varepsilon\in(0,1)$ such that $\psi_n>n^{\varepsilon}$, and hence $\varphi_n< n^{1-\varepsilon},$ so that that $\varphi_n$ is sublinear. 
In addition, the second inequality entails that $\psi_n$ is sublinear, too.
\tc{For the moment we assume that $B(\cdot)$ be non-lattice; see 
Remark~\ref{REM2} for the \tc{corresponding} result in the lattice case.}

The procedure we follow to prove the exact asymptotics in the slow regime is in line with the one we developed for the fast regime: we perform a change of measure, take out the exponential factor $\delta_n$, and analyze the remainder term $\Delta_n$. 
As this procedure echoes the one developed in the previous section, we only include the main steps and . 

\tc{The change of measure is again based on the solution $\vt_n$ that solves equation (\ref{twist}). 
In this case, as we argue below, $\vt_n$ obeys the expansion
\begin{equation}\label{theta_slow}
\vt_n = \sum_{k=1}^\infty w_k\psi_n^{\,-k};
\end{equation}
the coefficients $w_k$ can be recursively determined. The reasoning is analogous to the argumentation followed in the fast regime.}

\begin{itemize}
\item[$\circ$] \tc{To show the validity of the expansion (\ref{theta_slow}),
we write $\vt_-(x)$ as the solution of the equation
\[\beta'\left(\frac{\alpha(\vt)}{x}\right)\alpha'(\vt) =u,\]
so that $\vt_n= 
\vt_-(1/\psi_n)$ Recall that in this regime $\psi_n\to \infty $ as $n\to\infty$.}

\noindent \tc{
Observe that $\vt_-(x)/x\to \tau\s$ as $x\to 0$, where $\tau\s$ solves $a\beta'(a\tau\s)=u.$  We write $\vt_-(\cdot)$ as a Taylor expansion around $x=0$:
\[\vt_-(x) =\sum_{k=1}^\infty w_k x^k,\]
with $w_1=\tau\s$.
As in the fast regime, the coefficients allow expressions in terms of the derivatives of $\alpha(\cdot)$ and $\beta(\cdot)$.
Combining the above, we thus obtain the expansion (\ref{theta_slow}):
\[\vt_n=\vt_-(1/\psi_n)= \sum_{k=1}^\infty w_k \psi_n^{\,-k}.\]
}

\item[$\circ$] \tc{The procedure to identify the coefficients $w_k$ in this slow regime is analogous to the procedure to  identify the coefficients $v_k$ in this fast regime. In particular, 
$w_1$ solves
$a\beta'(aw_1) = u$;
in the sequel, we refer to $w_1$ by $\tau\s.$ }
\end{itemize}

We define $\sigma^{\mathbb Q}_-:=a\sqrt{\beta''(a\tau\s)}$;
cf.\ the decomposition (\ref{VAR}). 
For this regime we define
\[
\bar E_n := \frac{C_n - u n}{\psi_n\sqrt{\varphi_n}\sigma_-^{\mathbb Q}},
\]
which is a random variable that has, by construction, mean $0$ and variance converging to $1$ under ${\mathbb Q}_n$; comparing $\bar E_n$ with $\bar D_n$, observe that there is \tcp{a factor} $\psi_n\sqrt{\varphi_n}$ rather than $\sqrt{n}$ in the denominator. 
 Again we obtain, for all $n$, the factorization
\begin{equation}\label{ident_slow}
\xi_n(u) = {\rm e}^{\gamma_n(\vt_n)-\vt_n u n}\, \Delta_n,\:\:\:\mbox{with}\:\:\Delta_n:= {\mathbb E}_{{\mathbb Q}_n}\left( {\rm e}^{-\vt_n\sigma_-^{\mathbb Q}\psi_n\sqrt{\varphi_n}\,\bar E_n}\,1\{\bar E_n \geqslant 0 \}\right).\end{equation}
As before,  it remains to analyze  the exponential factor $\delta_n:=\exp(\gamma_n(\vt_n)-\vt_n u n)$ and the expectation $\Delta_n$, in the regime of $n$ growing large.

The analysis of $\delta_n$ precisely follows the corresponding step in the fast regime. 
Define 
\[
m_-:=\sup\left\{k\in{\mathbb N}: \liminf_{n\to\infty}\varphi_n\psi_n^{-k}>0\right\}.
\]
It thus follows that as $n\to\infty$, for constants ${\bar w}_k$, with the empty sum being defined as $0$,
\[
\delta_n=
\gamma_n(\vt_n)-\vt_n u n=
\big(\,\beta(a\tau\s)-\tau\s u\,\big)\varphi_n +\sum_{k=1}^{m_-} {\bar w}_k \varphi_n\psi_n^{-k}+o(1).
\]
\tcp{For the analysis of $\Delta_n$ we refer to Appendix \ref{Delta} and we proceed to Theorem \ref{THM_SLOW}. }

\begin{theorem}\label{THM_SLOW}
As $n\to\infty$, under Assumption \ref{ass2}, \tc{for non-lattice $B(\cdot)$,}
\[\xi_n(u)\sim \frac{1}{\tau\s \sigma_-^{\mathbb Q} \sqrt{2\pi\varphi_n}}
\exp\left(\big(\,  \beta(a\tau\s)-\tau\s u\,\big)\varphi_n +\sum_{k=1}^{m_-} {\bar w}_k \varphi_n\psi_n^{\, -k}\right).\]
\end{theorem}

\begin{corollary}\label{C2}
If  $\varphi_n\psi_n^{-1}=n/\psi_n^2  \to 0$ as $n\to\infty$, then, under Assumption \ref{ass2}, \tc{for non-lattice $B(\cdot)$,}
\[\xi_n(u) \sim {\mathbb P}(a B(\varphi_n)\geqslant u\varphi_n) \sim   \frac{1}{\tau\s \sigma_-^{\mathbb Q} \sqrt{2\pi\varphi_n}}
\exp\left(\big(\,  \beta(a\tau\s)-\tau\s u\,\big)\varphi_n \right).\]
As $n\to\infty$, under Assumption \ref{ass2}, \tc{for non-lattice $B(\cdot)$,}
\[\frac{1}{\varphi_n}\log\xi_n(u) \to  \beta(a\tau\s)-\tau\s u.\]
\end{corollary}

\begin{remark}\label{REM2}\tc{\em
In the case $B(\cdot)$ is lattice, the result of Thm.\ \ref{THM_SLOW} has to be slightly adjusted; cf.\  Remark~\ref{REM1}.
Let, for some $x_0$ and $d$ and any $t\geqslant 0$, the random variable $d^{-1}(B(t) -x_0)$  be integer almost surely, and let $d$ be the largest number with this property. Then, under Assumption \ref{ass2}, we obtain
\[\xi_n(u)\sim \frac{d}{1-{\rm e}^{-a\tau\s d}}\frac{1}{ \sigma_-^{\mathbb Q}/a\cdot \sqrt{2\pi\varphi_n}}
\exp\left(\big(\,  \beta(a\tau\s)-\tau\s u\,\big)\varphi_n +\sum_{k=1}^{m_-} {\bar w}_k \varphi_n\psi_n^{\, -k}\right),\]
as $n\to\infty$. Analogously with what we observed in the fast regime, this behavior is consistent with that of Thm.\ \ref{THM_SLOW} when taking $d\downarrow 0.$ $\hfill\Diamond$
}
\end{remark}

\begin{remark}\label{RemOverTekenA}
\tc{\em
In our analysis we have assumed that $a>0$; one may wonder whether our findings extend to general $a$. When picking $a\leqslant 0$, however, in the slow regime complications arise. The constant $\tau\s$, that plays a crucial role in this regime, should solve the equation $u=\beta'(a\tau\s)a$, but observe that $\beta'(\cdot)$ is positive (as it the characteristic exponent of a subordinator); as a consequence,  $u=\beta'(a\tau\s)a$ has no solution for $a\leqslant 0.$ This fact implies that the case $a\leqslant 0$ should be dealt with with an entirely different technique, in the sense that  the change-of-measure technique (as was used above) cannot work. 
$\hfill\Diamond$
}
\end{remark}

\section{Examples}\label{EX}

In this section we include two examples that demonstrate how the asymptotic expansion can be evaluated. 

\subsection{$A(\cdot)$ is a Poisson process and $B(\cdot)$ a Gamma process}
Let $A(\cdot)$ be a Poisson process; we assume its rate is $\lambda>0$, so that
$\alpha(\vt) = \lambda({\rm e}^{\vt} - 1)$ \tc{for $\vt\in{\mathbb R}$}. Let $B(\cdot)$ be a Gamma process; we call the parameters $r>0$ (shape) and $\mu>0$ (rate), so that $\beta(\vt) =r  \log\mu-r\log ({\mu - \vt})$, \tc{where we require $\vt<\mu$.}
Observe that  $A(t)$ has a Poisson distribution with parameter $\lambda t$, and that $B(t)$ has a Gamma distribution with parameters $rt$ and $\mu$. 
In particular, in the terminology of our paper, $a = \lambda$ and $b = r/\mu$. To make the event of interest {\it rare}, we assume $u>ab$; 
writing $\varrho :=\lambda r/({\mu u})$ this translates to assuming $\varrho <1$. \tc{Note that we have that $A(\cdot)$ is lattice but $B(\cdot)$ is not, so that we should apply Remark \ref{REM1} in the fast regime, and Thm.\ \ref{THM_SLOW} in the slow regime.}

We start our computations by providing the l-mgf of $A(\psi_n B( \phi_n))$ for this special case:
\begin{equation}\label{LLMGF}
\gamma_n(\vt)={\varphi_n } \,\beta\big( \alpha(\vt)\,\psi_n \big) = r\,{\varphi_n }\, \log\left(\frac{\mu}{\mu - \lambda({\rm e}^{\vt} - 1)\psi_n}\right).
\end{equation}
As $\vartheta_n$ satisfies the first-order condition $\gamma'_n(\vt_n) = u n$, we are to solve
\begin{equation}\label{thetaa}
\beta'(\alpha(\vt_n)\psi_n) \, \alpha'(\vt_n) = \frac{r\lambda\,{\rm e}^{\vt_n}}{\mu - \lambda({\rm e}^{\vt_n} - 1)\psi_n} = u.
\end{equation}
We thus find 
\[\vt_n = \log\left(\frac{\mu u+\lambda\psi_n u}{\lambda r+\lambda \psi_n u}\right) = \log\left(1 + \frac{\lambda}{\mu}\psi_n \right) - \log\left(\varrho \,( 1 +\frac{u}{r}\psi_n ) \right).\]
We \tcp{distinguish} between the fast regime and the slow regime. In the fast regime, in which $\psi_n\to 0$ as $n\to\infty$,  applying the Taylor expansion of the logarithm yields an expression for the coefficients $v_k$.
With $\zeta_1:=\lambda/\mu$ and $\zeta_2:=u/r$,
\[\vt_n=
\vt\s +\sum_{k=1}^\infty v_k \psi_n^{\,k},\:\:\:\:\vt\s = \log \frac{1}{\varrho},\:\:\:\:v_k:=
\frac{(-1)^{k+1}}{k} \left(\zeta_1^{\,k}-\zeta_2^{\,k}\right).
\]
We compute the coefficients $\bar v_k$ featuring in Lemma~\ref{EXP}. To this end, note that by inserting the first-order condition (\ref{thetaa}) into (\ref{LLMGF}), 
\[
\gamma_n(\vt_n) = -r \,{\varphi_n }\,\log(\varrho {\rm e}^{\vt_n})=- r\,{\varphi_n }\, (\vt_n - \vt\s) = -r \, \sum_{k=1}^{\infty} v_k \,\varphi_n\psi_n^{\,k},
\]
  so that
  \[\delta_n = \gamma_n(\vt_n) -\vt_n u n = -r\,\varphi_n \,\sum_{k=1}^{\infty} v_k \psi_n^k
  -\vt\s un - un \sum_{k=1}^\infty v_k\psi_n^{\,k}.\]
\tcp{It follows that} $- r \,v_1 
= ( 1- \varrho)u = b \alpha(\vt\s)$. 
 We thus observe that $\delta_n$ indeed has the form that was established in Lemma~\ref{EXP}, with, for $k \in\{  2,3,\ldots\}$, $\bar v_k=-(r\,v_k+u\,v_{k-1})$. More explicitly,
\[\bar{v}_k =
(-1)^{k} \left(\frac{r\big(\zeta_1^{\,k} -\zeta_2^{\,k} \big)}{k} - \frac{u \big(\zeta_1^{\,k-1} -\zeta_2^{\,k-1} \big)}{k-1}\right).\] 
Noting that $(\sigma_+^{\mathbb Q})^2 = b\alpha''(\vt\s) = u$, and bearing in mind Remark \ref{REM1},  we conclude that
\begin{equation}\label{xifast}\xi_n(u)\sim \frac{1}{1-\varrho}\frac{1}{ \sqrt{2\pi \,un }}
\exp\left(\left(\,
 1-\varrho
+ \log {\varrho}\right)u n +\sum_{k=2}^{m_+} {\bar v}_k \varphi_n\psi_n^{\,k}\right),\end{equation}
as $n\to\infty$.

The computations pertaining to the slow regime work very similarly. 
The crucial step is that we rewrite $\vartheta_n$ by `Tayloring' with respect to $\psi_n^{-1}$ (rather than to $\psi_n$): with $\bar \zeta_1:=\mu/\lambda$ and $\bar\zeta_2:=r/u$,
\[\vt_n = \log\left(1+\frac{\mu}{\lambda}\psi_n^{-1}\right)-  \log\left(1+\frac{r}{u}\psi_n^{-1}\right)=
\sum_{k=1}^\infty w_k \psi_n^{-k},\:\:\:w_k=\frac{(-1)^{k+1}}{k} \left(\bar\zeta_1^{\,k}-\bar\zeta_2^{\,k}\right).\]
Then $\tau\s = w_1 = \bar\zeta_1-\bar\zeta_2= \bar\zeta_1(1-\varrho) = \frac{r}{u}(\frac{1}{\varrho}-1)$ and $\beta(a \tau\s) = r \, \log\frac{1}{\varrho}$.
This gives
\begin{align*}
\gamma_n(\vt_n) &=
-r \,{\varphi_n } \log(\varrho\, {\rm e}^{\theta_n})= \beta(a \tau\s) \varphi_n - r\, \varphi_n \sum_{k=1}^\infty w_k \psi_n^{\,-k},
\end{align*}
so that 
 \[\delta_n = \gamma_n(\vt_n) -\vt_n u n = (\beta(a \tau\s) - \tau\s u) \varphi_n - \sum_{k=1}^\infty (r\,w_k + u\, w_{k+1}) \varphi_n \psi_n^{-k}.\]
\tc{Defining $(\sigma^{\mathbb{Q}}_-)^2 = a^2 \beta''(a \tau\s) ={u^2}/{r}$, application of Thm.\ \ref{THM_SLOW} leads after minor calculations to}
\begin{equation}
\label{xislow}\tc{\xi_n(u)\sim \frac{1}{\tfrac{1}{\varrho}-1}\frac{1}{ \sqrt{2\pi \,r \varphi_n }}
\exp\left(\left( 1- {\textstyle\frac{1}{\varrho}}+ \log{\textstyle\frac{1}{\varrho}}) \right)r\varphi_n +\sum_{k=1}^{m_-} {\bar w}_k \varphi_n\psi_n^{\,-k}\right),}\end{equation}
as $n\to\infty$, where $\bar{w}_k = -(r \,w_k +u\, w_{k+1})$ for $k \in\{ 1,2,\ldots\}$; observe the similarity with (\ref{xifast}).

\begin{remark}\tc{\em 
\label{REM3}
It is noted that the expressions encountered in this example align with those featuring in 
\cite[Section 3]{HKM} that deal with a case in which $\xi_n(u)$ could be evaluated explicitly.
It can be checked that $C_n$ has a negative binomial distribution, as follows. Observe that
\[\gamma_n(\vt) = r\varphi_n\log\left(\frac{\mu}{\mu-\lambda({\rm e}^{\vt}-1)\psi_n}\right).\]
In case a random variable  has a negative binomial distribution with parameters $k$ and $p$, then
its l-mgf is of the form 
\[
k \log\left(\frac{p}{1-(1-p)\mathrm{e}^{\vt}}\right).
\]
It thus follows that in the setting considered in this subsection, we have that $C_n$ is negative binomial, with success probability equal to $p= {\mu}/({\mu + \lambda \psi_n})$ and the allowed number of failures $k = r \varphi_n$.
Further intuition behind the appearance of the negative binomial distribution has been provided in
\cite[Remark 6]{HKM}.
\hfill$\Diamond$}
\end{remark}

\subsection{$A(\cdot)$ is a Gamma process and $B(\cdot)$ a Poisson process}
In our second  example the Gamma process and the Poisson process swap roles. 
\tc{In other words, we have $\alpha(\vt) =r  \log\mu-r\log ({\mu - \vt})$ for $\vt<\mu$, and $\beta(\vt) =
\lambda({\rm e}^{\vt} - 1)$ for $\vt\in{\mathbb R}$.} Again, to make the event of interest rare,  we   assume that $u>ab$, or alternatively, $\varrho := {\lambda r}/{(\mu u)} < 1$ (where $a=r/\mu$ and $b=\lambda$). \tc{In this case $A(\cdot)$ is non-lattice, whereas $B(\cdot)$ is, so that we have to use Thm.\ \ref{THM1} in the fast regime, and Remark \ref{REM2} in the slow regime.}
\begin{remark}{\em 
\label{REM4}
Note that the random variable $A(\psi_n B(\varphi_n))$ is compound Poisson in this case, with Gamma distributed jumps. More precisely, the jumps are generated according to a Poisson process with 
rate $\varphi_n \lambda$, where the jumps are Gamma, with parameters $r \psi_n$ and ${\mu}$.
\hfill$\Diamond$}
\end{remark}
In this case the l-mgf is given by
\begin{equation}
\label{gamma11}\gamma_n(\vt) =\varphi_n \lambda\left(\left(\frac{\mu}{\mu-\vt}\right)^{r\psi_n}-1\right).\end{equation}

Again we have to distinguish between the fast regime and the slow regime. 
As before, we start with the fast regime. 
From (\ref{gamma11}) it follows, by solving the first-order condition, that 
\[\vt_n = \mu\left(1-\varrho^{\eta_n}
\right),\:\:\:\eta_n:=
\frac{1}{1+r\psi_n},\]
so that
\[\gamma_n(\vt_n) = \varphi_n\lambda\left(\textstyle\frac{1}{\varrho} \cdot \varrho^{\eta_n}-1\right).\]

To find the $v_k$, we write $\vt_n$ as a Taylor series in $\psi_n$.
It is directly seen
that $\vt\s = \mu(1 - \varrho);$
some elementary calculus thus yields for
the first coefficients $v_1=
-\mu r\varrho\, \log\frac{1}{\varrho}$, and \[v_2= \mu r^2 \varrho\,(\log \textstyle\frac{1}{\varrho})\,(1 - \frac{1}{2}\log \textstyle\frac{1}{\varrho});\]
higher coefficients can be found in the same way.

Also the coefficients $\bar v_k$ can be found: 
\begin{align*}
\delta_n &= \gamma_n(\vt_n) - \vt_n u n = \varphi_n \lambda \left(\textstyle \frac{1}{\varrho} \cdot \varrho^{\eta_n} - 1\right) - \mu (1 - \varrho^{\eta_n}) u n \\
&=\frac{\lambda}{\varrho}(\varrho^{\eta_n}-1+1-\varrho)\varphi_n-\mu(1-\varrho^{\eta_n})un=\frac{\lambda}{\varrho}(1-\varrho)\varphi_n -\mu u \left(\textstyle\frac{1}{r}\varphi_n+n\right)(1-\varrho^{\eta_n})\\
&=\frac{\lambda}{\varrho}(1-\varrho)\varphi_n -u \left(\textstyle\frac{1}{r}\varphi_n+n\right)
\left(\vt\s+\sum_{k=1}^\infty v_k \psi_n^k\right)\\
&=\left(\frac{\lambda}{\varrho}-\frac{\mu u}{r}\right) (1-\varrho) \varphi_n - (1-\varrho)\mu u n -{\frac{u}{r}v_1 n} - u
\sum_{k=2}^\infty\left(\textstyle\frac{1}{r}v_k +v_{k-1}\right)\varphi_n\psi_n^{\,k}\\
&= \left(1 - \textstyle\frac{1}{\varrho} + \log\textstyle\frac{1}{\varrho}\right) \lambda r n - u\sum_{k=2}^{\infty} \left(\textstyle\frac{1}{r} v_k +   v_{k-1}\right) \varphi_n \psi_n^k, 
\end{align*}
so that $\bar{v}_k = -u\,( \frac{1}{r} v_k +  v_{k-1})$.
Note that the first term equals $(b \alpha(\vt\s) - \vt\s u ) n$, in accordance with the result in Lemma~\ref{EXP}.
\tc{Thm.\ \ref{THM1}, with $(\sigma_+^{\mathbb Q})^2 = u^2/(\lambda r)$, yields}
\[\xi_n(u)\sim \frac{1}{(\frac{1}{\varrho} - 1) }\frac{1}{\sqrt{2\pi  \,\lambda r n}}
\exp\left(\left(1 - \textstyle\frac{1}{\varrho} + \log\textstyle\frac{1}{\varrho}\right)\lambda r n +\sum_{k=2}^{m_+} {\bar v}_k \varphi_n\psi_n^{\,k}\right).\]

We proceed with the slow case. 
The $w_k$ follow by expanding $\vt_n$ in terms of a Taylor series  with respect to $\psi_n^{-1}$ (recalling that $\psi_n\to\infty$):
\[
\vt_n= 
\mu\left(1- \varrho \cdot \varrho^{\bar\eta_n} \right),\:\:\:\:\bar\eta_n:=-\frac{r}{r+1/{\psi_n}}.
\]
Routine calculations show that $\tau\s = w_1 = ( {\mu}/{r})\,\log \frac{1}{\varrho} $, and 
\[w_2 = -\frac{\mu}{r^2}\, (\log\textstyle\frac{1}{\varrho})\,(1 + \frac12 \log\textstyle \frac{1}{\varrho}),\]
where higher coefficients follow along the same lines.
The $\bar{w}_k$ can be determined as before; leaving out a few intermediate steps,
\begin{align*}
\delta_n &= \gamma_n(\vt_n) - \vt_n u n = \varphi_n \lambda (\varrho^{\bar\eta_n} -1)- \mu (1 - \varrho\cdot \varrho^{\bar\eta_n}) u n \\
&=  \left (1-\varrho + \log\varrho\right) \frac{\lambda}{\varrho}\,\varphi_n -u \sum_{k=1}^{\infty} (\textstyle\frac{1}{r}w_k + w_{k+1}) \varphi_n \psi_n^{\,-k}.
\end{align*}
\tc{By applying Remark \ref{REM2}, with $(\sigma_-^{\mathbb Q})^2 = (r/\mu)^2 \cdot \lambda/\varrho = ru/\mu$ and $\bar{w}_k = - u\,(\frac{1}{r} w_k +  w_{k+1})$,
\[\xi_n(u)\sim \frac{1}{1- \tfrac{1}{\varrho}}\frac{1}{ \sqrt{2\pi \,(\lambda / \varrho)\, \varphi_n}}
\exp\left(\big(\, 1- \varrho + \log {\varrho}\,\big)\,({\lambda}/{\varrho})  \,\varphi_n +\sum_{k=1}^{m_-} {\bar w}_k \varphi_n\psi_n^{\,-k}\right).\]}

\section{Applications in modeling overdispersed customer streams}
\label{OVERD}
In operations research, the performance evaluation and design of service systems is a key topic of interest.
Traditionally, the dominant assumption when modeling customer arrival processes is that of Poisson arrivals. 
As we mentioned in the introduction, however, measurements indicate that the level of variation observed in practice may significantly exceed that predicted by the Poisson process  \cite{BAS,LIU,MAT, WHI}. More concretely, where for Poisson arrivals the mean and variance of the number of arrivals in a given time interval coincide, in practice one typically observes that the variance is larger than the mean; the phenomenon is typically referred to as {\it overdispersion}.  

In \cite{HLM} a mechanism is proposed that produces overdispersion. The main idea is that the arrival rate is random rather than deterministic; this is achieved by {\it resampling} the arrival rate periodically. The main underlying idea behind the resampling is that there is a (perhaps unobservable) environmental process that influences all potential customers  in the same way; think for instance of the weather conditions, or the occurrence of specific events.  As we argued in the introduction, 
our two-timescale framework can be used to represent this `resampled Poisson process'. 

Renormalizing time such that the resampling intervals have length 1, considering $K$ such intervals, and letting the samples be i.i.d.\ non-negative random variable $\Lambda_1,\ldots,\Lambda_K$, then the number of arrivals in there $K$  intervals is distributed as 
\[C(K) = {\rm Pois}\left(\sum_{i=1}^K \Lambda_i\right).\]
When designing service centers (for instance when making staffing decisions), one would like to control the probability of excessive waiting times. For that reason, it could be informational to have a handle on the tail distribution corresponding to the number of arrivals in a certain time window. 

The remainder of this section has two objectives. In the first place, we quantify the error made by neglecting overdispersion (i.e., the situation in which the number of arrivals in the $K$ slots has a Poisson distribution with parameter $K\,{\mathbb E}\Lambda$). In the second place, we compare our refined asymptotics with the crude asymptotics that were found in \cite{HLM}, just involving the dominant term in the exponent (as featuring in the second statements of Corollaries \ref{C1} and \ref{C2}).

\subsection{Setup}
In our experiments, the following example is studied. 
Let $\Lambda_1,\ldots,\Lambda_K$ be i.i.d.\ samples from 
an exponential distribution with mean $1/\bar\mu$. We are interested in using our expansions to approximate the probability
\begin{equation}\Pi(K,\bar u,\bar\mu):={\mathbb P}\left(C(K) \geqslant \bar u\right),\label{Pint}\end{equation}
for some $\bar u>K/\bar\mu.$ 
\begin{itemize}
\item[$\circ$]
Importantly, as we have seen in Remark \ref{REM3}, this example allows an explicit solution: here $C(K)$ has a negative binomial distribution with parameters $K$ and success probability $\bar\mu/(1+\bar\mu)$.  We thus have, in self-evident notation,
\[\Pi(K,\bar u,\bar\mu) = {\mathbb P}\left({\rm NegBin}\left(K,\frac{\bar\mu}{1+\bar\mu}\right)\geqslant \bar u\right).\]
We chose the example of exponentially distributed $\Lambda_i$, as the explicit expressions allow us to compare the various approximations with the corresponding exact values. 
\item[$\circ$]
When neglecting the overdispersion, one does not work with sampled values $\Lambda_1,\ldots,\Lambda_K$, but rather with a deterministic arrival rate ${\mathbb E}\Lambda$ in any interval. In other words, when one would ignore the overdispersion, the probability $\Pi(K,\bar u,\bar\mu)$ would be approximated by
\[\Pi^{\rm (Pois)}(K,\bar u,\bar\mu):= {\mathbb P}\left({\rm Pois}\left(K\,{\mathbb E} \Lambda\right)\geqslant \bar u\right).\]
This approximation, corresponding with the fast regime, is expected to work well if $K$ is large relative to $\bar u$, and the contributions of the individual $\Lambda_i$ are relatively small (i.e., $\bar\mu$ is relatively large).
\item[$\circ$]\tc{
Along the same lines, one could consider the regime that $K$ is substantial, but not in the order of $\bar u$, and the contributions of the individual $\Lambda_i$ are relatively large (i.e., $\bar\mu$ is relatively small). \tcp{Here} the slow regime is expected to apply, suggesting to ignore the Poisson sampling. Bearing in mind that the sum of exponentially distributed random variables (with equal parameter) has an Erlang distribution, we end up with (in self-evident notation)
\[\Pi^{\rm (Gamma)}(K,\bar u,\bar\mu):= {\mathbb P}\left({\rm Erl}\left(K,\bar\mu\right)\geqslant \bar u\right).\]}
\end{itemize}

The crude approximations $\Pi^{\rm (Pois)}(K,\bar u,\bar\mu)$ and $\Pi^{\rm (Gamma)}(K,\bar u,\bar\mu)$ we wish to compare with approximations based on our two-timescale model. To this end, we convert the asymptotics that we found in Section \ref{EX} into approximations. 
We fix a number $\tr{n>1}$; this scaling parameter can be considered as `artificial' in the sense that the  choice of $n$ will not affect the approximation (as will turn out below). Then we pick $f$ such that $n^{f}=K$ (i.e., $f= \log K/\log n>0$), and we set $u := \bar u/n$ and $ \mu:=\bar\mu n^{1-f}$.  We
thus obtain
\[\sum_{i=1}^K \Lambda_i\stackrel{\rm d}{=} n^{1-f} \sum_{i=1}^{n^{f}} \Lambda_i^\circ,\]
where $\Lambda_i^\circ$ is exponentially distributed with parameter $\mu.$
We have thus rewritten the probability (\ref{Pint}) as 
\[{\mathbb P}\left({\rm Pois}\left(n^{1-f} \sum_{i=1}^{n^{f}}  \Lambda_i^\circ\right)\geqslant un\right).\]
\tcp{Notice} that this representation falls in the framework of Section 4.1, with $\varphi_n = n^f$. 
More specifically, we are in the situation that $A(\cdot)$ is a Poisson process with rate $\lambda =1$, and $B(\cdot)$ a Gamma process with shape parameter $r=1$ and rate parameter $\mu$.

The next step is to translate the asymptotics that we identified for $\xi_n(u)$, which are in terms of the parameters  of the scaled model (i.e., $\mu$, $f$, and $u$), into approximations for $\Pi(K,\bar u,\bar\mu)$, which are in terms of the original parameters (i.e., $\bar\mu$, $K$, and $\bar u$). 

\subsection{Fast regime}
We start by considering the regime in which $K$ is large relative to $\bar u$, so that it makes sense to work with the asymptotics pertaining to the case $f>1$ (i.e., the fast regime). 
First observe that
\begin{equation}
\label{rootje}\varrho=\frac{\lambda r}{\mu u} = \frac{K}{\bar\mu \,\bar u},\end{equation}
which we assumed to be smaller than 1. In addition,
\[\zeta_1 = \frac{1}{\bar \mu \,n^{1-f}} = \frac{1}{\bar\mu\psi_n},\:\:\:\:\: \zeta_2 = \frac{\bar u}{n} = \frac{\bar u}{K\psi_n},\]
so that (after some elementary computations)
\begin{equation}
\label{teetje}\bar v_k \varphi_n \psi_n^{\,k} = (-1)^k \left( \frac{K\,t^+_k}{k} 
-\frac{\bar u\,t^+_{k-1}}{k-1}\right),\:\:\:\:\:
t^+_k:= \left(\frac{1}{\bar\mu}\right)^k - \left(\frac{\bar u}{K}\right)^k.\end{equation}
Also,
\[\vt\s = \log\frac{\bar\mu\,\bar u}{K}.\]
Based on (\ref{xifast}), we eventually obtain the following approximation for the probability (\ref{Pint}):
\begin{equation}
\label{APPR1}\Pi^{\rm (fast)}(K,\bar u,\bar\mu):=\frac{1}{1-\varrho}\frac{1}{\sqrt{2\pi \bar u}}
\exp\left((1-\varrho+\log\varrho)\, \bar u + \sum_{k=2}^{M_+} (-1)^k \left( \frac{K\,t^+_k}{k} 
-\frac{\bar u\,t^+_{k-1}}{k-1}\right)\right),\end{equation}
with $\varrho$ given by (\ref{rootje}) and $t_k^+$ by (\ref{teetje}).
In practice one should choose the threshold $M_+$ such that adding more  terms in the sum  does not affect the outcome; realize that the object $m_+$ is not well-defined in the pre-limit context. Observe that, as announced, the scaling parameter $n$ does not appear in the approximation. 

In \cite{HLM} logarithmic asymptotics of $\xi_n(u)$ were found. These lead to (crude) approximations that 
only cover the dominant term in the exponent (as in the second statement of Corollary \ref{C1}). More concretely, one would approximate the probability of our interest by
\[\hat\Pi^{\rm (fast)}(K,\bar u,\bar\mu)=
\exp\left((1-\varrho+\log\varrho)\, \bar u \right).\]

\subsection{Slow regime}
Now consider the regime in which $K$ is small relative to $\bar u$, in which we opt for relying on the asymptotics corresponding to $f<1.$ Then 
\[\bar\zeta_1 = \bar\mu\,n^{1-f} = \bar\mu \,\psi_n,\:\:\:\:\bar\zeta_2= \frac{n}{\bar u} =\frac{K\psi_n}{\bar u},\]
implying that
\begin{equation}\label{teetjemin}\bar w_k\varphi_n \psi_n^{\,-k} =(-1)^k\left( \frac{K\,t^-_k}{k} 
-\frac{\bar u\,t^+_{k+1}}{k+1}\right),\:\:\:\:\:
t^-_k:= {\bar\mu}^k - \left(\frac{K}{\bar u}\right)^k.
\end{equation}
In addition, 
\[\tau\s = {K^{1/f}}\left(\frac{\bar\mu}{K}-\frac{1}{\bar u}\right).\]
Based on (\ref{xislow}), we thus obtain the following approximation for (\ref{Pint}):
\begin{equation}
\label{APPR2}\Pi^{\rm (slow)}(K,\bar u,\bar\mu):= \frac{1}{\tfrac{1}{\varrho}-1}\frac{1}{\sqrt{2\pi K}}\exp\left(\left(1-\tfrac{1}{\varrho}+\log\tfrac{1}{\varrho}\right)K+\sum_{k=1}^{M_-}(-1)^k\left( \frac{K\,t^-_k}{k} 
-\frac{\bar u\,t^-_{k+1}}{k+1}\right)
\right),\end{equation}
with $\varrho$ given by (\ref{rootje}) and $t_k^-$ by (\ref{teetjemin}). As before, the truncation level $M_-$ is chosen such that including additional terms would not have any impact on the approximation. Again, we observe that the scaling parameter $n$ does not appear in the approximation.

We also compare with the approximation that 
only covers the dominant term in the exponent (as in the second statement of Corollary \ref{C2}, in line with the decay rate identified in \cite{HLM}). 
This approximation reads
\[\hat\Pi^{\rm (slow)}(K,\bar u,\bar\mu)=\exp\left(\left(1-\tfrac{1}{\varrho}+\log\tfrac{1}{\varrho}\right)K\right).\]

\subsection{Numerical experiments}
We conclude this section by presenting a set of examples that illustrate the use of the  developed approximations. We subsequently consider an example that corresponds to the fast regime, and one that corresponds to the slow regime. They provide indications of the accuracy that can be achieved by the expressions (\ref{APPR1}) and (\ref{APPR2}) that were based on our refined asymptotics, relative to more crude approximations.

\vspace{3mm}

We start with an example corresponding to the fast regime. Recalling (\ref{APPR1}), notice that our approximation $\Pi^{\rm (fast)}(K,\bar u,\bar\mu)$
is parameterized by the threshold $M_+$. In Table \ref{T1} below, the column with superscript 0 corresponds to neglecting the sum in the exponent of (\ref{APPR1}) altogether, which could be seen as choosing $M_+:=1$, whereas the column with superscript 1 corresponds to including one term in the sum in the exponent of (\ref{APPR1}), i.e.,  $M_+:=2$. Throughout this example  $\varrho = \tfrac{1}{3}$ and $\bar u= 150$ are held fixed, explaining that the values in the third, fourth and fifth column are constant (as can be seen from the expressions for the approximations $\Pi^{\rm (Pois)}, \hat\Pi^{\rm (fast)}$, and $\Pi^{\rm (fast, 0)}$). 

\begin{table}[h]{\normalsize 
\begin{tabular}{|| c |ccccc||}\hline
{$(K, \bar\mu,\bar u)$} & $\Pi$ & $\Pi^{\rm (Pois)}$ & $\hat\Pi^{\rm (fast)}$ & $\Pi^{\rm (fast, 0)}$ & $\Pi^{\rm (fast, 1)}$ \\
\hline
$ (100\,000, 1\,000, 150)$& $1.90\cdot10^{-6}$ &$1.88\cdot 10^{-6}$  & $2.00\cdot 10^{-5}$ &$1.95\cdot10^{-6}$& $1.97\cdot10^{-6}$ \\
$ (50\,000, 500, 150)$& $1.93\cdot10^{-6}$ &$1.88\cdot 10^{-6}$  & $2.00\cdot 10^{-5}$  &$1.95\cdot10^{-6}$& $2.00\cdot10^{-6}$ \\
$ (10\,000, 100, 150)$ &$2.13\cdot10^{-6}$ & $1.88\cdot 10^{-6}$ & $2.00\cdot 10^{-5}$ &$1.95\cdot10^{-6}$ &$2.21\cdot10^{-6}$ \\
$ (5\,000, 50, 150)$& $2.41\cdot10^{-6}$ &$1.88\cdot 10^{-6}$  & $2.00\cdot 10^{-5}$  &$1.95\cdot10^{-6}$& $2.51\cdot10^{-6}$ \\
$ (1\,000, 10, 150)$ & $5.89\cdot10^{-6}$ & $1.88\cdot 10^{-6}$ & $2.00\cdot 10^{-5}$ &$1.95\cdot10^{-6}$ & $6.82\cdot10^{-6}$\\ \hline
\end{tabular}\vspace{2mm}
}
\caption{Approximations for fast regime.\label{T1}}
\end{table}
The following observations can be made  from Table \ref{T1}. (i)~In the \tr{two} top rows there is so much timescale separation that $\Pi^{\rm (fast, 0)}$ provides highly accurate results. In the \tr{other} rows the timescale separation becomes less pronounced, which leads to $\Pi^{\rm (fast, 1)}$ becoming the preferred approximation. Here it is noted that including 2 or more terms in the exponent hardly improves the approximation in the third and fourth row. In the last row the timescales are so poorly separated that $\Pi^{\rm (fast, 1)}$ is still relatively far off; we mention that the performance is significantly improved by using $M_+=3$, which leads to the highly accurate  approximation $5.90\cdot10^{-6}.$ (ii)~As anticipated, the Poisson approximation (neglecting the overdispersion) performs well if there is a substantial degree of timescale separation, but leads to substantial underestimation if the timescales are relatively close together. 

\vspace{3mm}

\tr{In the slow regime  our approximation $\Pi^{\rm (fast)}(K,\bar u,\bar\mu)$, as given by  (\ref{APPR2}),
is parameterized by the threshold $M_-$. Analogously to the notation used in Table \ref{T1}, in Table \ref{T2} the column with superscript 0 corresponds to neglecting the sum in the exponent of (\ref{APPR2}), which one could identify with the choice $M_-:=0$, whereas the column with superscript 1 corresponds to including one term in the sum in the exponent of (\ref{APPR2}), i.e.,  $M_-:=1$. 
We chose scenarios in which $\varrho=\tfrac{2}{3}$ and $K=100$ are held fixed, implying that the values in the fourth and fifth column are constant (as can be seen from the expressions for $\hat\Pi^{\rm (fast)}$ and $\Pi^{\rm (fast, 0)}$); note that in this slow regime the values in the third column (i.e., $\Pi^{\rm (Gamma)}$) are not constant, in that they (very slowly)  increase.
}

%
%
%

\begin{table}[h!]{\normalsize 
\begin{tabular}{|| c |ccccc||}\hline
{$(K, \bar\mu,\bar u)$} & $\Pi$ & $\Pi^{\rm (Gamma)}$ & $\hat\Pi^{\rm (fast)}$ & $\Pi^{\rm (fast, 0)}$ & $\Pi^{\rm (fast, 1)}$ \\
\hline
$(100, 0.001,150\,000) $& $5,98\cdot 10^{-6}$ & $5,93\cdot 10^{-6}$ & $7,84\cdot 10^{-5}$ &$6,26\cdot 10^{-6}$& $6,31\cdot 10^{-6}$ \\
$ (100, 0.05,30\,000)$ & $6.20\cdot 10^{-6}$ & $5.93\cdot 10^{-6}$& $7.84\cdot 10^{-5}$  & $6.26\cdot 10^{-6}$ & $6.52\cdot 10^{-6}$ \\
$ (100, 0.01,15\,000)$ &  $6.48\cdot 10^{-6}$& $5.95\cdot 10^{-6}$ &  $7.84\cdot 10^{-5}$ & $6,26\cdot 10^{-6}$ &$6.80\cdot 10^{-6}$ \\
$(100, 0.05,3000) $ & $9.11\cdot 10^{-6}$ & $6.03\cdot 10^{-6}$ & $7.84\cdot 10^{-5}$ & $6,26\cdot 10^{-6}$ & $9.49\cdot 10^{-6}$\\
$(100, 0.1,1500) $ & $1.35\cdot 10^{-5}$ & $6.14\cdot 10^{-6}$ & $7.84\cdot 10^{-5}$ & $6,26\cdot 10^{-6}$ & $1.44\cdot 10^{-5}$\\
\hline
\end{tabular}\vspace{2mm}
}
\caption{Approximations for slow regime.\label{T2}}
\end{table}

\tr{The conclusions from  Table \ref{T2} are as follows. (i)~Again, in the top rows (with a strong timescale separation) $\Pi^{\rm (fast, 0)}$ performs well. This approximation, however, degrades in the lower rows, where the accuracy substantially improves if $M_-=1$ (i.e., the approximation $\Pi^{\rm (fast, 1)}$) is used. We mention that for the parameters in the last row, the accuracy is further improved by taking $M_-=2$, yielding $1.34\cdot 10^{-5}.$ (ii)~The Gamma approximation (assuming full separation of time scales) provides accurate approximations in the top row, but substantially worse in the bottom rows (where there is less timescale separation), as expected.}

\section{Discussion and concluding remarks}\label{DISC}
\tr{
Motivated by recent developments in arrival process modeling, this paper presents tail asymptotics for a multi-timescale model. The focus has been on approximating $\xi_n(u)={\mathbb P}(A(\psi_nB(\varphi_n)\geqslant un)$,
with $A(\cdot)$ a L\'evy process and $B(\cdot)$ a L\'evy subordinator. Our analysis shows that subtle analysis is required to identify the exact asymptotics. The structure found is considerably richer than in the model's classical (single timescale) counterpart; in particular, the exponent includes additional terms.}

\subsection*{Effect of leaving customers}
\tr{In the case $A(\cdot)$ is a Poisson process, we recover the mixed Poisson model proposed in \cite{HLM}, which can be thought of as a Poisson process with periodically resampled rate (so as to generate overdispersion). An interesting extension of the results developed in the presented paper could concern the model in which the arrived clients leave after a service time distributed as the random variable $J$ with $F_J(t):={\mathbb P}(J\leqslant t)$. Such a setting could be considered as an infinite-server queue with overdispersed input. Using the insights from \cite{HLM}, one thus obtains the following expression for the number of customers $\bar C_n$ present at time $\varphi_n$, under the scaling considered in the present paper:
\[\bar C_n= A\left(\psi_n\int_0^{\varphi_n} (1-F_J(t))\,{\rm d} B(t)\right).\]
Using calculation rules for L\'evy processes, one obtains
\[\log {\mathbb E}{\rm e}^{\vt \bar C_n} = \int_0^{\varphi_n} \beta\big(\alpha(\vt) \psi_n (1-F_J(t))\big) {\rm d}t;\]
plugging in $F_J(t)=0$ for all $t\in[0,\varphi_n]$, we recover ${\varphi_n} \beta(\alpha(\vt) \psi_n)$ (which makes sense, as for this choice customers do not leave the system). As a topic for further research, one could pursue deriving the exact asymptotics of ${\mathbb P}(\bar C_n\geqslant un)$ (which may require scaling the service times $J$). These asymptotics can be converted into approximations, to be used as the basis of refined staffing rules.}

\subsection*{Follow-up research}
\tr{Other directions for future work include: 
\begin{itemize}
\item[$\circ$] In the first place, multivariate extensions could be considered. One could for instance consider the probability that the vector
\[\big(A_1(\psi_nB(\varphi_n)), A_2(\psi_nB(\varphi_n))\big),\]
with ${\boldsymbol A}(\cdot)=(A_1(\cdot),A_2(\cdot))$ a bivariate L\'evy processes, attains a value in the set $[u_1n,\infty)\times [u_2n,\infty)$. The  components of ${\boldsymbol A}(\cdot)$ could be assumed dependent, but observe that  $A_1(\psi_nB(\varphi_n))$ and $A_2(\psi_nB(\varphi_n))$ are dependent anyway (as the same $B(\varphi_n)$ is used). The multivariate single-timescale results of \cite{CS} show that one should expect that various cases arise: there will be cases in which one component exceeding its threshold (with high probability) implies that the other component exceeds its threshold as well, but also cases in which both constraints are `tight'. 
\item[$\circ$] In Remark \ref{RemOverTekenA} we pointed out that in our setup it is required to assume $a>0$, particularly in the slow regime. For $a\leqslant 0$, one may anticipate that $\xi_n(u)$ decays exponentially in $n$ (unlike in the case of $a>0$, where $\xi_n(u)$ decays exponentially in $\varphi_n$), based on the following reasoning. The starting point is the large-deviations heuristic
\begin{equation}\label{LDH}\xi_n(u) \approx \max_{x>0} {\mathbb P}(B(\varphi_n) \approx x \varphi_n) \,{\mathbb P}(A(xn) \approx un).\end{equation}
The main insight is that, taking into account only the leading term of the asymptotics, the first probability on the right-hand side of (\ref{LDH}) decays exponentially in $\varphi_n$, but the second exponentially in $n$, so that their product decays exponentially in $n$. This should be contrasted with  the case $a>0$ studied in the present paper: there the maximum will be attained (roughly) at $x=u/a>0$, so that $\{A(xn) \approx un\}$ is no rare event, and hence ${\mathbb P}(aB(\varphi_n) \approx u\varphi_n)$ dictates the tail behavior (thus leading to exponential decay in $\varphi_n$). A comparable situation has been considered in \cite[Section 2.2]{HKM}.
We stress however that, relative to the $a>0$ case, this $a\leqslant 0$ case has a considerably lower practical relevance. 
\end{itemize}}

\bibliographystyle{plain}

\appendix

\section{Edgeworth expansions} \label{Edge}
In this appendix we establish Edgeworth expansions for $C_n$. 
We successively address the fast and the slow regime.
\subsection{Fast regime} 
The goal is to develop an expansion for ${\mathbb Q}_n(\bar D_n\leqslant x)$. The proof is a variation of that for sums of i.i.d.\ random variables \cite{FEL}, and therefore we only provide the main steps. 
First observe that
\begin{equation}
\label{expa}{\mathbb E}_{{\mathbb Q}_n}\,{\rm e}^{\vt\bar D_n} = \big({\rm e}^{\Gamma_n(\vt/\sigma_+^{\mathbb Q}\sqrt{n})}\big)^{\varphi_n},\end{equation}
with $\Gamma_n(\vt):= \beta(\alpha(\vt+\vt_n)\psi_n) -\beta(\alpha(\vt_n)\psi_n)-\vt u\psi_n.$ Obviously, 
\[{\rm e}^{\Gamma_n(\vt/(\sigma_+^{\mathbb Q}\sqrt{n}))} = \sum_{k=0}^\infty \frac{1}{k!} \omega_n^{(k)},\:\:\:\:
 \omega_n^{(k)}:=\left(\Gamma_n
 \left(\frac{\vt}{\sigma_+^{\mathbb Q}\sqrt{n}}\right)\right)^k.\]
Due to the definition of $\vt_n$,  $\Gamma_n'(0)=0$. It requires some elementary calculus to verify that, with \tr{$\tr{\gamma^\circ_+}:=\beta''(0)(\alpha'(\vt\s))^2 +\beta''(0)\alpha(\vt\s)\alpha''(\vt\s)+\beta'(0)\alpha'''(\vt\s)v_1$},
\begin{align*}
\Gamma_n''(0) &= \beta''(\alpha(\vt_n)\psi_n)\big(\alpha'(\vt_n)\big)^2\psi_n^2+\beta'\big(\alpha(\vt_n)\psi_n\big)\alpha''(\vt_n)\psi_n\\
&= b \alpha''(\vt\s)\psi_n +\tr{\tr{\gamma^\circ_+}}\psi_n^2+o(\psi_n^2),
\end{align*}
and $\Gamma_n'''(0) = b\alpha'''(\vt\s)\psi_n+o(\psi_n)$.
Upon combining the above results, we obtain
\begin{align*}
{\rm e}^{\Gamma_n(\vt/(\sigma_+^{\mathbb Q}\sqrt{n}))} =\:& 1+\frac{\vt^2}{2\varphi_n}+
\frac{1}{2}\frac{\tr{\gamma^\circ_+}}{(\sigma_+^{\mathbb Q})^2}\frac{\vt^2\psi_n}{\varphi_n}\,+\\
&\frac{1}{6}\frac{b\alpha'''(\vt\s)}{(\sigma_+^{\mathbb Q})^3}\frac{\vt^3}{\varphi_n\sqrt{n}}+o\left(\max\left\{\frac{\psi_n}{\varphi_n},\frac{1}{\varphi_n\sqrt{n}}\right\}\right).
\end{align*}
Note that we have to distinguish between two cases: $\lim_{n \to \infty }\psi_n \sqrt{n}= 0$ and $\liminf_{n \to \infty} \psi_n \sqrt{n} > 0$.

In the case where $\lim_{n \to \infty} \psi_n \sqrt{n} = 0$, we have due to (\ref{expa}),
\[{\mathbb E}_{{\mathbb Q}_n}\,{\rm e}^{\vt\bar D_n} =
\left(1+ \frac{\vt^2}{2\varphi_n}+
\frac{1}{6}\frac{b\alpha'''(\vt\s)}{(\sigma_+^{\mathbb Q})^3}
\frac{\vt^3}{\varphi_n\sqrt{n}}+o\left(\frac{1}{\varphi_n\sqrt{n}}\right)\right)^{\varphi_n},\]
which can be rewritten as
\[\left(1+\frac{\vt^2}{2\varphi_n}\right)^{\varphi_n} +
\left(1+ \frac{\vt^2}{2\varphi_n}\right)^{\varphi_n-1} 
\frac{1}{6}\frac{b\alpha'''(\vt\s)}{(\sigma_+^{\mathbb Q})^3}\frac{\vt^3}{\sqrt{n}}+o\left(\frac{1}{\sqrt{n}}\right).\]
\tc{By a direct computation}, we find 
\[
{\mathbb E}_{{\mathbb Q}_n}\,{\rm e}^{\vt\bar D_n} =
\exp\left(\frac{1}{2}{\vt^2}\right)\left(1+
\frac{1}{6}\frac{b\alpha'''(\vt\s)}{(\sigma_+^{\mathbb Q})^3}
\frac{\vt^3}{\sqrt{n}} \right)+o\left(\frac{1}{\sqrt{n}}\right).
\]
Using the \tc{familiar} inversion procedure for characteristic functions, we thus obtain, with $\phi(\cdot)$ denoting the probability density function of a standard Normal random variable,\[{\mathbb Q}_n(\bar D_n\in {\rm d}x) = \phi(x)\left(1+H_3(x)\frac{1}{6}\frac{b\alpha'''(\vt\s)}{(\sigma_+^{\mathbb Q})^3}\frac{1}{\sqrt{n}}\right)+
o\left(\frac{1}{\sqrt{n}}\right),\]
with $H_k(\cdot)$ the Hermite polynomial of degree $k$. 
This leads to, with $\Phi(\cdot)$ denoting the cumulative distribution function of a standard Normal random variable,
\begin{equation} \label{bound1fast}
{\mathbb Q}_n(\bar D_n\leqslant x) = \Phi(x)-\phi(x)\left(H_2(x)\frac{1}{6}\frac{b\alpha'''(\vt\s)}{(\sigma_+^{\mathbb Q})^3}\frac{1}{\sqrt{n}}\right)+
o\left(\frac{1}{\sqrt{n}}\right).
\end{equation}
With the same reasoning as in the standard Edgeworth expansion (i.e., that for sums of i.i.d.\ random variables), the error term (being small relative to $1/\sqrt{n}$) is uniform in $x$. 

In the case where $\liminf_{n \to \infty} \psi_n \sqrt{n}>0$, by inserting our expansion into (\ref{expa}),
\[
{\mathbb E}_{{\mathbb Q}_n}\,{\rm e}^{\vt\bar D_n} =
\left(1+\frac{\vt^2}{2\varphi_n}+
\frac{1}{2}\frac{\tr{\gamma^\circ_+}}{(\sigma_+^{\mathbb Q})^2}\frac{\vt^2\psi_n}{\varphi_n}+o\left(\frac{\psi_n}{\varphi_n}\right)\right)^{\varphi_n},
\]
which can be rewritten as
\[
\left(1+\frac{\vt^2}{2\varphi_n}\right)^{\varphi_n} +
\left(1+\frac{\vt^2}{2\varphi_n}\right)^{\varphi_n-1} 
\frac{1}{2}\frac{\tr{\gamma^\circ_+}}{(\sigma_+^{\mathbb Q})^2}{\vt^2\psi_n}+o\left(\psi_n\right).
\]
Using the same procedure as above, this yields, uniformly in $x$,
\[
{\mathbb Q}_n(\bar D_n\leqslant x) = \Phi(x)-\phi(x)\left(H_1(x)\frac{1}{2}\frac{\tr{\gamma^\circ_+}}{(\sigma_+^{\mathbb Q})^2}\psi_n \right)+
o\left(\psi_n\right).
\]
In our setting, however, we need an error that is small relative to $1/\sqrt{n}$, which can be achieved by expanding $\Gamma_n''(0)$ further (as it can be verified  that the contributions of the higher derivatives $\Gamma_n^{(k)}(0)$ for $k\geqslant 3$ are $o(1/\sqrt{n})$). Define
\[
k_+ := \sup\left\{k\in{\mathbb N}: \liminf_{n\to\infty} {\psi_n^{\, k}}{\sqrt{n}}>0\right\};
\]
due to Assumption \ref{ass1}, this is a finite constant.
Using the above type of reasoning, we conclude that there are constants $c_1,\ldots,c_{k_+}$ so that
\begin{equation} \label{bound2fast}
{\mathbb Q}_n(\bar D_n\leqslant x) = \Phi(x)-\phi(x)\left(H_1(x)\sum_{k=1}^{k_+}c_k\psi_n^{\, k} +H_2(x)\frac{1}{6}\frac{b\alpha'''(\vt\s)}{(\sigma_+^{\mathbb Q})^3}\frac{1}{\sqrt{n}}\right)+
o\left(\frac{1}{\sqrt{n}}\right).
\end{equation}
{In the boundary case where $\psi_n \sqrt{n}$ converges to a constant we have $k_+=1$; the sum in (\ref{bound2fast}) consists of only one term.}

\subsection{Slow regime} 
\tcp{Here} our objective is to find an expansion for ${\mathbb Q}_n(\bar E_n\leqslant x)$; the reasoning is analogous to that of the fast regime. 
The starting point is the mgf of $\bar E_n$ under the twisted distribution:
\begin{equation} \label{expa2}
{\mathbb E}_{{\mathbb Q}_n}\,{\rm e}^{\vt\bar E_n} = \big({\rm e}^{\Gamma_n(\vt/\sigma_-^{\mathbb Q}\sqrt{n\psi_n})}\big)^{\varphi_n},
\end{equation}
with, as before, $\Gamma_n(\vt):= \beta(\alpha(\vt+\vt_n)\psi_n) -\beta(\alpha(\vt_n)\psi_n)-\vt u\psi_n$. 
\tcp{We get}
\begin{equation}
\label{ps}\frac{1}{2} \left(\frac{\vt}{\sigma_-^{\mathbb Q}\sqrt{n\psi_n}}\right)^2 \Gamma_n''(0)=
\frac{1}{2} \frac{\vt^2}{(\sigma_-^{\mathbb Q})^2} \left(
\beta''(\alpha(\vt_n)\psi_n)\big(\alpha'(\vt_n)\big)^2\frac{1}{\varphi_n}+\beta'\big(\alpha(\vt_n)\psi_n\big)\alpha''(\vt_n)\frac{1}{n}\right).\end{equation}
\tr{Using the expansion of $\vt_n$, the right-hand side of the previous display reads
\[
\frac{1}{2} \frac{\vt^2}{(\sigma_-^{\mathbb Q})^2} \left(
\beta''(a\tau\s)a^2\frac{1}{\varphi_n}+\frac{\gamma^\circ_-}{n}+o\left(\frac{1}{n}\right)\right)=\frac{\vt^2}{2\varphi_n}+\frac{1}{2}\frac{\gamma^\circ_-}{(\sigma_-^{\mathbb Q})^2}\frac{\vt^2}{n}+o\left(\frac{1}{n}\right),
\]
where $\gamma^\circ_-:=\beta'\big(a\tau\s\big)\alpha''(0)+2a\, \alpha''(0)\beta''(a\tau\s)\tau\s+\tfrac{1}{2}a^2\alpha''(0)  \beta'''(a\tau\s)\,(\tau\s)^2
+a^3\beta'''(a\tau\s)w_2.$
}In addition,
\[
\frac{1}{6} \left(\frac{\vt}{\sigma_-^{\mathbb Q}\sqrt{n\psi_n}}\right)^3 \Gamma_n'''(0)=
\frac{1}{6}\frac{\beta'''(a\tau\s)a^3}{(\sigma_-^{\mathbb Q})^3}\frac{\vt^3}{\varphi_n^{3/2}}+o\left(\frac{1}{\varphi_n^{3/2}}\right).
\]
We distinguish between the cases $\lim_{n \to \infty} \varphi_n^{3/2}/n = 0$ and $\liminf_{n \to \infty} \varphi_n^{3/2}/n > 0$. 
Mimicking the reasoning used for the fast regime, in the case where $\lim_{n \to \infty} \varphi_n^{3/2}/n = 0$,
\[{\mathbb E}_{{\mathbb Q}_n}\,{\rm e}^{\vt\bar E_n} =
\left(1+ \frac{\vt^2}{2\varphi_n}+\frac{1}{6}\frac{\beta'''(a\tau\s)a^3}{(\sigma_-^{\mathbb Q})^3}\frac{\vt^3}{\varphi_n^{3/2}}+o\left(\frac{1}{\varphi_n^{3/2}}\right)\right)^{\varphi_n},\]
which can be rewritten as
\[\left(1+\frac{\vt^2}{2\varphi_n}\right)^{\varphi_n} +
\left(1+ \frac{\vt^2}{2\varphi_n}\right)^{\varphi_n-1} 
\frac{1}{6}\frac{\beta'''(a\tau\s)a^3}{(\sigma_-^{\mathbb Q})^3}\frac{\vt^3}{\sqrt{\varphi_n}}+o\left(\frac{1}{\sqrt{\varphi_n}}\right).\]
This leads to, inserting our expansion into (\ref{expa2}),
\[{\mathbb E}_{{\mathbb Q}_n}\,{\rm e}^{\vt\bar E_n} =
\exp\left(\frac{1}{2}{\vt^2}\right)\left(1+
\frac{1}{6}\frac{\beta'''(a\tau\s)a^3}{(\sigma_-^{\mathbb Q})^3}\frac{\vt^3}{\sqrt{\varphi_n}}
\right)+o\left(\frac{1}{\sqrt{\varphi_n}}\right),\]
and, after inversion, to the Edgeworth expansion
\begin{equation} \label{bound1slow}
{\mathbb Q}_n(\bar E_n\leqslant x) = \Phi(x)-\phi(x)\left(H_2(x)\frac{1}{6}\frac{\beta'''(a\tau\s)a^3}{(\sigma_-^{\mathbb Q})^3}\frac{1}{\sqrt{\varphi_n}}\right)+
o\left(\frac{1}{\sqrt{\varphi_n}}\right).
\end{equation}
Finally, focus on $\liminf_{n \to \infty}\varphi_n^{3/2}/n > 0$. 
Performing the same steps,
\[
{\mathbb E}_{{\mathbb Q}_n}\,{\rm e}^{\vt\bar E_n} =
\exp\left(\frac{1}{2}{\vt^2}\right)\left(1+
\frac{1}{2}\frac{\tr{\gamma^\circ_-}}{(\sigma_-^{\mathbb Q})^2}\frac{\vt^2}{\psi_n}+o\left(\frac{1}{\psi_n}\right)\right),
\]
leading to, uniformly in $x$,
\[
{\mathbb Q}_n(\bar E_n\leqslant x) = \Phi(x)-\phi(x)\left(H_1(x)\frac{1}{2}\frac{\tr{\gamma^\circ_-}}{(\sigma_-^{\mathbb Q})^2}\frac{1}{\psi_n} \right)+
o\left(\frac{1}{\psi_n}\right).
\]
Our objective, however, is to obtain an error that is $o(1/\sqrt{\varphi_n})$. This is achieved by expanding (\ref{ps}) further; again the contributions of  the higher derivatives (i.e., derivatives of order 3 and higher) can be verified to be negligible. Following the line of reasoning of the fast regime, we define
\[
k_- := \sup\left\{k\in{\mathbb N}: \liminf_{n\to\infty} \frac{{\sqrt{\varphi_n}}}{\psi_n^{\, k}}>0\right\};
\]
due to Assumption \ref{ass2}, this is a finite constant.
Using the above type of reasoning, we conclude that there are constants ${c'_1,\ldots,c'_{k_-}}$ so that 
\begin{equation}\label{bound2slow}
{\mathbb Q}_n(\bar E_n\leqslant x) = \Phi(x)-\phi(x)\left(H_1(x)\sum_{k=1}^{k_-}\frac{c'_k}{\psi_n^{\, k}} +H_2(x)\frac{1}{6}\frac{\beta'''(a\tau\s)a^3}{(\sigma_-^{\mathbb Q})^3}\frac{1}{\sqrt{\varphi_n}}\right)+
o\left(\frac{1}{\sqrt{\varphi_n}}\right).
\end{equation}
Again, in the boundary case where $\varphi_n^{3/2}/n$ converges to a constant we have $k_-=1$ and the sum in (\ref{bound2slow}) consists of only one term.

\section{Analysis of $\Delta$\textsubscript{\lowercase{n}} in the slow regime} \label{Delta} 
\tcp{Analogously to our analysis in the fast regime, applying integration by parts, we write $\Delta_n$ as}
\begin{align} \nonumber
\Delta_n =&\,\int_0^\infty {\rm e}^{-\vt_n\sigma_-^{\mathbb Q}\psi_n\sqrt{\varphi_n}\,x} {\mathbb Q}_n(\bar E_n\in{\rm d}x)\\ \nonumber =&\,
\psi_n\sqrt{\varphi_n}\vt_n\sigma_-^{\mathbb Q}\int_0^\infty {\rm e}^{-\vt_n\sigma_-^{\mathbb Q}\psi_n\sqrt{\varphi_n}\,x}\big( {\mathbb Q}_n(\bar E_n\le x)- {\mathbb Q}_n(\bar E_n\le 0)\big){\rm d}x
\\
=&\,\psi_n\vt_n\sigma_-^{\mathbb Q}\int_0^\infty {\rm e}^{-\psi_n\vt_n\sigma_-^{\mathbb Q} x} \big({\mathbb Q}_n(\bar E_n\leqslant x/\sqrt{\varphi_n})-
{\mathbb Q}_n(\bar E_n\leqslant 0)\big){\rm d}x.\label{ABC}
\end{align}
We again proceed by using the Edgeworth expansion presented in \tcp{Appendix \ref{Edge}}.

$\circ$\:\:In the case $\lim_{n \to \infty}\varphi_n^{3/2}/n = 0$ (which corresponds to $f<\frac23$ in the context of Example \ref{exf}), we have, as pointed out in Eqn.\ (\ref{bound1slow}) in \tcp{Appendix \ref{Edge}}, as $n\to\infty$, 
\begin{equation} \label{B3}
\sqrt{\varphi_n}\sup_x\left({\mathbb Q}_n(\bar E_n\leqslant x)- \Phi(x)+\phi(x) H_2(x)\frac{\kappa_-}{\sqrt{\varphi_n}}\right)\to 0,\:\:\:\kappa_-:=\frac{1}{6}\frac{ \beta'''(a \tau\s)a^3}{ (\sigma_-^{\mathbb Q})^3}.
\end{equation}
Our objective is to prove that $\sqrt{\varphi_n}\Delta_n$ converges to a constant; as in the fast regime, we provide the upper bound, but the lower bound follows fully analogously. 
By   (\ref{ABC}) and (\ref{B3}), for any $\varepsilon>0$ and $n$ sufficiently large, 
\begin{align*}
\sqrt{\varphi_n}\Delta_n \leqslant &\:\sqrt{\varphi_n}\psi_n\vt_n \sigma_-^{\mathbb Q} \int_0^\infty {\rm e}^{-\psi_n\vt_n\sigma_-^{\mathbb Q} x} \left(\Phi(\frac{x}{\sqrt{\varphi_n}})-\Phi(0)\right)
{\rm d}x\,-\\
&\:\sqrt{\varphi_n}\psi_n\vt_n\sigma_-^{\mathbb Q}  \int_0^\infty {\rm e}^{-\psi_n\vt_n\sigma_-^{\mathbb Q} x} \left(\phi(\frac{x}{\sqrt{\varphi_n}})H_2(\frac{x}{\sqrt{\varphi_n}})-\phi(0)H_2(0)\right)\frac{\kappa_-}{\sqrt{\varphi_n}}
{\rm d}x\,+\,\varepsilon.
\end{align*}
\tc{The first term on the right-hand side equals}
\[ 
{\sqrt{\varphi_n}}\exp\left(\frac{1}{2}(\vt_n\sigma_-^{\mathbb Q} \psi_n )^2 \varphi_n\right) \left(1-\Phi(\psi_n\vt_n\sigma_-^{\mathbb Q} \sqrt{\varphi_n})\right),
\]
which, using $x(1-\Phi(x))/\phi(x)\to 1$ and $\psi_n\vt_n\to\tau\s$,  converges to
\begin{equation} \label{constanteslow}
\big({\tau\s \sigma_-^{\mathbb Q} \sqrt{2\pi}}\big)^{-1}.
\end{equation}
As before, we split the second term  into (recalling that $H_2(x)=x^2-1$)
\begin{align*}
t_{-,n}^{(1)}&:=\,\sqrt{\varphi_n}\psi_n\vt_n \sigma_-^{\mathbb Q} \int_0^\infty {\rm e}^{-\psi_n\vt_n\sigma_-^{\mathbb Q} x} \left(\phi(0)-\phi(\frac{x}{\sqrt{\varphi_n}})\right)\frac{\kappa_-}{\sqrt{\varphi_n}}
{\rm d}x{,}\\
t_{-,n}^{(2)}&:=\,
\sqrt{\varphi_n}\psi_n\vt_n \sigma_-^{\mathbb Q} \int_0^\infty {\rm e}^{-\psi_n\vt_n\sigma_-^{\mathbb Q} x} \phi(\frac{x}{\sqrt{\varphi_n}})\frac{x^2}{\varphi_n}\frac{\kappa_-}{\sqrt{\varphi_n}}
{\rm d}x.
\end{align*}
Mimicking the reasoning used in the fast regime, it can be shown that $t_{-,n}^{(1)}\to 0$ and
 $t_{-,n}^{(2)}\to 0$ as $n\to\infty$.
The third term equals $\varepsilon$, which can be made arbitrarily small. 
Combining this with  the corresponding lower bound applies, we find that $\sqrt{\varphi_n}\Delta_n$ converges to  (\ref{constanteslow}).

$\circ$\:\: {Considering} $\liminf_{n \to \infty} \varphi_n^{3/2}/n > 0$  (which simplifies to $f\in [\frac{2}{3},1)$ in the context of Example~\ref{exf}), we obtain from Eqn.\ (\ref{bound2slow}) in \tcp{Appendix \ref{Edge}},
\begin{equation} \label{B4}
\sqrt{\varphi_n}\sup_x\left({\mathbb Q}_n(\bar E_n\leqslant x)- \Phi(x)+\phi(x) \left(H_1(x)\sum_{k=1}^{k_-}c'_k\psi_n^{-k} +H_2(x)\frac{\kappa_-}{\sqrt{\varphi_n}}\right)\right)\to 0.
\end{equation}
By (\ref{B4}), for any $\varepsilon>0$ and $n$ sufficiently large, using that $H_1(x)=x$ and hence $H_1(0)=0$,
\begin{align*}
\sqrt{\varphi_n}\Delta_n \leqslant &\:\sqrt{\varphi_n}\psi_n\vt_n \sigma_-^{\mathbb Q} \int_0^\infty {\rm e}^{-\psi_n\vt_n\sigma_-^{\mathbb Q} x} \left(\Phi(\frac{x}{\sqrt{\varphi_n}})-\Phi(0)\right)
{\rm d}x\,-\\
&\:\sqrt{\varphi_n}\psi_n\vt_n \sigma_-^{\mathbb Q} \int_0^\infty {\rm e}^{-\psi_n\vt_n\sigma_-^{\mathbb Q} x}  \phi(\frac{x}{\sqrt{\varphi_n}})\frac{x}{\sqrt{\varphi_n}}\left(\sum_{k=1}^{k_-}c_k'\psi_n^{-k} \right) 
{\rm d}x\,-\,\\
&\:\sqrt{\varphi_n}\psi_n\vt_n \sigma_-^{\mathbb Q} \int_0^\infty {\rm e}^{-\psi_n\vt_n\sigma_-^{\mathbb Q} x} \left(\phi(\frac{x}{\sqrt{\varphi_n}})H_2(\frac{x}{\sqrt{\varphi_n}})-\phi(0)H_2(0)\right)\frac{\kappa_-}{\sqrt{\varphi_n}}
{\rm d}x\,+\,\varepsilon.
\end{align*}
The first, third, and fourth term can be dealt with as in the case $\varphi_n^{3/2}/n\to  0$. 
Analogously to the reasoning used in the fast regime, the second term vanishes. The corresponding lower bound is established in the same way. 
We conclude that also in this case $\sqrt{\varphi_n}\Delta_n$ converges to  (\ref{constanteslow}). 

\end{document}